\documentclass[leqno]{amsart}
\usepackage{amssymb}
\usepackage{amsbsy}
\usepackage{amscd}
\usepackage{amsmath}
\usepackage{amsthm}
\usepackage[all]{xy}

\newtheorem{theo}{Theorem}[section]
\newtheorem{defi}{Definition}\numberwithin{defi}{section}

\newtheorem{assu}{Assumption}

\newtheorem{lemm}{Lemma}
\newtheorem{prop}{Proposition}

\newtheorem{rema}{Remark}
\newtheorem{coro}{Corollary}

\newtheorem{exa}[theo]{Example}

\newcommand{\sogot}{\ensuremath{\mathfrak{so}}}

\newcommand{\Acal}{\ensuremath{\mathcal{A}}}

\newcommand{\Dcal}{\ensuremath{\mathcal{D}}}
\newcommand{\Ecal}{\ensuremath{\mathcal{E}}}
\newcommand{\Fcal}{\ensuremath{\mathcal{F}}}

\newcommand{\Hcal}{\ensuremath{\mathcal{H}}}
\newcommand{\Kcal}{\ensuremath{\mathcal{K}}}

\newcommand{\Pcal}{\ensuremath{\mathcal{P}}}
\newcommand{\Qcal}{\ensuremath{\mathcal{Q}}}
\newcommand{\Scal}{\ensuremath{\mathcal{S}}}
\newcommand{\Vcal}{\ensuremath{\mathcal{V}}}

\newcommand{\Cbb}{\ensuremath{\mathbb{C}}}
\newcommand{\Rbb}{\ensuremath{\mathbb{R}}}
\newcommand{\Zbb}{\ensuremath{\mathbb{Z}}}

\newcommand{\Kbb}{\ensuremath{\mathbb{K}}}

\newcommand{\A}{\ensuremath{\mathbb{A}}}
\newcommand{\F}{\ensuremath{\mathbf{F}}}

\newcommand{\Ahat}{\ensuremath{\widehat{A}}}

\newcommand{\dr}{\ensuremath{\hbox{\scriptsize \rm dec-rap}}}
\newcommand{\KK}{\ensuremath{\mathbf K}}
\newcommand{\sm}{\ensuremath{\rm sm}}
\newcommand{\x}{\ensuremath{\mathbf{x}}}
\newcommand{\cbf}{\ensuremath{\mathbf{c}}}

\newcommand{\f}{\ensuremath{\mathcal{C}^{\infty}}}

\newcommand{\T}{\ensuremath{\hbox{\bf T}}}
\newcommand{\End}{\ensuremath{\hbox{\rm End}}}
\newcommand{\str}{\operatorname{Str}}
\newcommand{\herm}{\operatorname{Herm}}
\newcommand{\tr}{\operatorname{Tr}}
\newcommand{\ch}{\operatorname{Ch}}
\newcommand{\chc}{\operatorname{Ch_{\rm c}}}
\newcommand{\Kch}{\operatorname{\bf Ch}}

\newcommand{\chg}{\operatorname{Ch_{\rm sup}}}

\newcommand{\chq}{\operatorname{Ch_{\rm Q}}}

\newcommand{\expo}{\operatorname{exp}}
\newcommand{\chr}{\operatorname{Ch_{\rm rel}}}
\newcommand{\supp}{\operatorname{Supp}}
\newcommand{\p}{\operatorname{p}}

\newcommand{\Eul}{\operatorname{Eul}}
\newcommand{\Pf}{\operatorname{Pf}}
\newcommand{\bere}{\operatorname{T}}
\newcommand{\e}{\operatorname{e}}

\def \cf {\hbox{\scriptsize \rm fiber cpt}}

\def \id {{\rm Id}}

\def \tuc {{\rm Th_c}}

\def \tur {{\rm Th_{{\ensuremath{\rm rel}}}}}

\def \tumq {{\rm Th_{{\ensuremath{\rm MQ}}}}}

\title{Quillen's relative Chern character is multiplicative}

\author{Paul-Emile Paradan}

\address{Institut de Math\'ematiques et Mod\'elisation de Montpellier (I3M), 
Universit\'e Montpellier 2, 34095 Montpellier} 

\email{Paul-Emile.Paradan@math.univ-montp2.fr}

\author{Mich\`ele Vergne}

\address{Institut de Math\'ematiques de Jussieu, Th{\'e}orie des
Groupes, 75005 Paris -- Ecole Polytechnique,  Centre de math\'ematiques Laurent Schwartz,
91128 Palaiseau}

\email{vergne@math.polytechnique.fr}

\date{September 2008}

\begin{document}

\begin{abstract}
In the first part of this paper we prove the multiplicative property of the  
relative Quillen Chern character. Then we obtain a Riemann-Roch formula between 
the relative Chern character of the Bott morphism and
the relative Thom form.
\end{abstract}

\maketitle

{\small \tableofcontents}

\section{Introduction}

The relative Chern character was defined by Atiyah and al. in
\cite{Atiyah-Hirz61,Atiyah-Singer-2}
as a map
\begin{equation} \label{eq:ch-intro1}
\Kch_{X\backslash Y} :\KK^0(X,Y)\longrightarrow {\rm H}^*(X,Y).
\end{equation}
Here $Y\subset X$ are {\em finite} CW-complexes, $\KK^0(X,Y)$ is the
relative K-group and ${\rm H}^*(X,Y)$ is the singular relative
cohomology group.

The relative Chern character enjoys various functorial properties.
In particular,  $\Kch$ is multiplicative: the following diagram
\begin{equation}\label{eq:fonctoriel-intro}
\xymatrix@C=8mm{ \KK^0(X,Y)\ar[d]^{\Kch_{X\backslash Y}} \!\!\!\!&
\times \quad \KK^0(X,Y')\ar[d]^{\Kch_{X\backslash Y'}}\quad
\ar[r]^{\odot} &
\quad \KK^0(X,Y\cup Y') \ar[d]_{\Kch_{X\backslash Y\cup Y'}}\\
{\rm H}^*(X,Y)  \!\!\!\!\!   & \times \quad {\rm H}^*(X,Y')\quad
\ar[r]^{\diamond} & \quad {\rm H}^*(X,Y\cup Y') }
\end{equation}
is commutative. Here $Y,Y'\subset X$ are {\em finite}
CW-complexes  and $\odot$ and $\diamond$ denote the products. This
property was extended to the case where $X$ is a paracompact
topological space and $Y$ any \emph{open} subset of $X$ by Iversen
in \cite{Iversen76} (see also
\cite{Illusie-Sem19-Cartan,Illusie-LNM}). Iversen deduces the
existence of the local  Chern character  from functorial
properties, but his construction is not explicit.

In this article, we work in the context of manifolds and
differential forms. Indeed, in this framework, Quillen
\cite{Quillen85} constructed a very natural de Rham relative
cohomology class associated to a smooth morphism between vector
bundles, that we call the relative Quillen Chern character. Let
$N$ be a manifold, and let $\sigma:\Ecal^+\to\Ecal^-$ be a
morphism of complex vector bundles over $N$. Let $\supp(\sigma)$
be the \emph{support} of $\sigma$ : it is the set of points $n\in
N$ where  $\sigma(n)$ is not invertible. We do not suppose that
$\supp(\sigma)$ is compact.  Quillen \cite{Quillen85} associates
to $\sigma$ a couple $(\alpha,\beta)$ of differential forms, where
$\alpha$ is given by the usual Chern-Weil construction, and
 $\beta$ is also constructed {\`a} la Chern-Weil, via super-connections.
The form  $\alpha$ is a closed differential form on $N$
representing the difference of  Chern characters
$\ch(\Ecal^+)-\ch(\Ecal^-)\in \Hcal^*(N)$, and $\beta$ is a
differential form on $N\setminus\supp(\sigma)$ such that
$$
\alpha|_{N\setminus\supp(\sigma)}=d\beta.
$$

The couple $(\alpha,\beta)$ defines then an explicit relative de
Rham cohomology class
$$
\chr(\sigma)\in \Hcal^*(N,N\setminus\supp(\sigma)).
$$

The main purpose of this note is to show that  Quillen's relative
Chern character $\chr$ is multiplicative.
 If $\sigma_1,\sigma_2$ are two morphisms on
$N$, then the product $\sigma_1\odot\sigma_2$ is a morphism on $N$
with support equal to $\supp(\sigma_1)\cap\supp(\sigma_2)$. We
prove in Section \ref{sec:chr-multiplicative} that the following
equality
\begin{equation}\label{eq:chr-intro}
    \chr(\sigma_1\odot\sigma_2)=\chr(\sigma_1)\diamond\chr(\sigma_2)
\end{equation}
holds in $\Hcal^*(N,N\setminus\supp(\sigma_1\odot\sigma_2))$.

\medskip

Intuitively (and  true in many analytic cases), the relative Chern
class could also be represented as a current supported on
$\supp(\sigma)$, but currents do not usually multiply. Thus
another procedure, involving   a choice of partition of unity, is
needed  to define the product $\diamond$ of relative classes in de
Rham relative cohomology. The multiplicativity property
(\ref{eq:chr-intro}) can also be deduced from the fact that
Quillen's Chern character gives an explicit representative of
Iversen's local Chern character, due to Schneiders functorial
characterization of Iversen's class (see \cite{Schneiders}). Our
proof does not use Iversen's construction, and our explicit
argument  can be extended to the case of equivariant Chern
characters  with generalized coefficients (see
\cite{pep-vergne3}).

When $\supp(\sigma)$ is \emph{compact}, there is a natural
homomorphism  from  $\Hcal^*(N,N\setminus \supp(\sigma))$ into the
compactly supported cohomology algebra $\Hcal_c^*(N)$ and the image
of the Quillen relative Chern character $\chr(\sigma)$ is the Chern
character $\ch_c(\sigma)$ with compact support. The equality
(\ref{eq:chr-intro}) implies the relation
\begin{equation}\label{eq:chc-intro}
\ch_c(\sigma_1\odot\sigma_2)=\ch_c(\sigma_1)\wedge\ch_c(\sigma_2)\quad
\mathrm{in}\quad \Hcal^*_c(N).
\end{equation}
This last relation is well known and follows also from the fact
that $\ch_c(\sigma)$ is the Chern character of a difference bundle
on a compactification of $N$.

%

As an important example, we consider $\sigma_b$  the Bott morphism
on a complex  vector bundle $\sigma_b:\Lambda^+\Vcal\to
\Lambda^-\Vcal$ over $\Vcal$, given by the exterior product by
$v\in \Vcal$. This morphism has support the zero section $M$ of
$\Vcal$. It leads to a   relative class $\chr(\sigma_b)$ in
$\Hcal^*(\Vcal,\Vcal\setminus M)$. One can  give a  similar
construction of the relative  Thom form $\tur(\Vcal)\in
\Hcal^*(\Vcal,\Vcal\setminus M)$
 of the vector bundle $\Vcal\to M$, using the Berezin integral instead of a
super-trace.   The explicit formulae  for  $\chr(\sigma_b)$ and
$\tur(\Vcal)$ allows us to derive the ``Riemann-Roch" relation
between these two relative classes at the level of differential
forms. Our  proof  follows the same scheme than the proof of the
relation  between  the Chern character and the Thom class with
Gaussian looks constructed by Mathai-Quillen
\cite{Mathai-Quillen}.

\bigskip

\textbf{Acknowledgements:} We are grateful to M. Karoubi,  J.
Lannes, P. Schapira and J.P. Schneiders for enlightening
discussions on these topics.

We wish to thank the referee for  his careful reading, and suggestions for improvements.

\section{Cohomological structures}

Let $N$ be a manifold. We denote by   $\Acal^*(N)$ the algebra of
differential forms on $N$ and by  $\Hcal^*(N)$ the de Rham
cohomology algebra of $N$. We denote by
$\Hcal^*_c(N)$ its compactly supported cohomology algebra.

In this paper, we work with differential forms with {\em complex} or {\em real} coefficients, 
depending on the context. In order to simplify the notation, we use the same notation for 
$\Acal^*$, $\Hcal^*$ and $\Hcal^*_c$ viewed as {\em complex} or {\em real} vector space :  
we speak of $\Kbb$-differential forms, $\Kbb$-cohomology classes  or $\Kbb$-algebras  with
$$\Kbb\in\{\Rbb,\Cbb\}.$$

\subsection{Relative cohomology}\label{sec:cohomologie-relative}

Let $F$ be a closed subset of $N$. To a cohomology class on $N$
vanishing on $N\setminus F$, we associate a relative cohomology
class. Let us explain the construction (see \cite{Bott-Tu}).
Consider the $\Kbb$-complex $\Acal^*(N,N\setminus F)$ with
$$
\Acal^k(N,N\setminus F):=\Acal^k(N)\oplus \Acal^{k-1}(N\setminus F)
$$
and differential $d_{\rm rel}\left(\alpha,\beta\right)=
\left(d\alpha,\alpha|_{N\setminus F}- d\beta \right)$.

\begin{defi}\label{relcoh}
The cohomology of the complex $(\Acal^*(N,N\setminus F),d_{\rm
rel})$ is the relative $\Kbb$-cohomology space $\Hcal^*(N,N\setminus F)$.
\end{defi}

The class defined by a $d_{\rm rel}$-closed element
$(\alpha,\beta)\in\Acal^*(N,N\setminus F)$ will be denoted by
$[\alpha,\beta]$.
There is a natural $\Kbb$-linear map $\Hcal^*(N,N\setminus F)\to \Hcal^*(N)$.

If $F_1$ and $F_2$ are closed subsets of $N$,  there is a
natural product
\begin{eqnarray}\label{eq:produit-relatif}
\Hcal^*(N,N\setminus F_1)\times\Hcal^*(N,N\setminus F_2)
&\longrightarrow&
\Hcal^*(N,N\setminus (F_1\cap F_2)) \\
\quad\quad\quad\quad(\quad a \quad,\quad b\quad
)\quad\quad\quad\quad &\longmapsto& a\diamond b \ ,\nonumber
\end{eqnarray}
which is $\Kbb$-bilinear.


We will use an explicit formula for $\diamond$ that we recall.
Let $U_1:=N\setminus
F_1$, $U_2:=N\setminus F_2$ so that $U:=N\setminus (F_1\cap
F_2)=U_1\cup U_2$. Let $\Phi:=(\Phi_1,\Phi_2)$ be a partition of
unity subordinate to the covering $U_1\cup U_2$ of $U$. With the help of
$\Phi$, we define a bilinear map
$\diamond_\Phi:\Acal^*(N,N\setminus F_1)\times
\Acal^*(N,N\setminus F_2)$ $\to \Acal^*(N,N\setminus (F_1\cap
F_2))$ as follows. For $a_i:=(\alpha_i,\beta_i)\in \Acal^{k_i}(N,N\setminus F_i)$, $i=1,2$, we
define
$$
a_1 \diamond_\Phi a_2 :=\Big(\alpha_1\wedge \alpha_2,
\Phi_1\beta_1 \wedge \alpha_2+(-1)^{k_1}\alpha_1\wedge
\Phi_2\beta_2-(-1)^{k_1}d\Phi_1\wedge \beta_1\wedge \beta_2\Big).
$$

Remark that all forms $\Phi_1\beta_1 \wedge \alpha_2$,
$\alpha_1\wedge \Phi_2\beta_2$ and $d\Phi_1\wedge \beta_1\wedge
\beta_2$ are well defined on $U_1\cup U_2$. Indeed the support of
the form $d\Phi_1$ is contained in $U_1\cap U_2$, as
$d\Phi_1=-d\Phi_2$. So $a_1 \diamond_\Phi a_2\in
\Acal^{k_1+k_2}(N,N\setminus (F_1\cap F_2))$.
It is immediate to verify that $d_{\rm rel}(a_1 \diamond_\Phi a_2)$ is equal to
$(d_{\rm rel}a_1) \diamond_\Phi a_2+(-1)^{k_1}a_1 \diamond_\Phi (d_{\rm rel}a_2)$.
Thus $\diamond_\Phi$ defines a bilinear map $\Hcal^*(N,N\setminus F_1)\times
\Hcal^*(N,N\setminus F_2)$ $\to \Hcal^*(N,N\setminus (F_1\cap F_2))$.

Let us see that this product do not depend on the choice of the
partition of unity. If we have another partition
$\Phi'=(\Phi_1',\Phi_2')$, then
$\Phi_1-\Phi'_1=-(\Phi_2-\Phi'_2)$. It is immediate to verify
that, if $d_{\rm rel}(a_1)=0$ and $d_{\rm rel}(a_2)=0$, one has
$$
a_1 \diamond_\Phi a_2 - a_1 \diamond_{\Phi'}a_2=d_{\rm
rel}\Big(0,(-1)^{k_1} (\Phi_1-\Phi'_1)\beta_1\wedge \beta_2\Big).
$$
So the product on the relative cohomology spaces will be denoted by
$\diamond$.

\subsection{Inverse limit of cohomology with support}\label{sec:cohomologie-support}

Let $F$ be a closed subset of $N$. We consider the set $\Fcal_F$ of
all open neighborhoods $U$ of $F$ which is ordered by the relation
$U\leq V$ if and only if $V\subset U$. For any $U\in \Fcal_F$, we
consider the $\Kbb$-algebra $\Acal^*_U(N)$ of differential forms on $N$
with support contained in $U$ (that is vanishing on a neighborhood
of $N\setminus U$): this algebra is stable under the de Rham
differential $d$, and we denote by $\Hcal^*_U(N)$ the corresponding
cohomology $\Kbb$-algebra. If $U\leq V$, we have then an inclusion map
$\Acal^*_V(N)\hookrightarrow \Acal^*_U(N)$ which gives rise to a $\Kbb$-linear map
$f_{U,V}:\Hcal^*_V(N)\to \Hcal^*_U(N)$.

\begin{defi}\label{defiinductive}
 We denote by   $\Hcal^*_F(N)$ the inverse limit of the inverse system
 $(\Hcal^*_U(N),f_{U,V};U,V\in\Fcal_F)$. It is a $\Kbb$-vector space.
\end{defi}

%

If $F_1, F_2$ are two closed subsets of $N$,  there is a
$\Kbb$-bilinear map
\begin{equation}\label{eq:produit-support}
    \Hcal^*_{F_1}(N)\times \Hcal^*_{F_2}(N)
    \stackrel{\wedge}{\longrightarrow} \Hcal^*_{F_1\cap F_2}(N)
\end{equation}
which is defined via the wedge product on forms.

\bigskip

Now we define a $\Kbb$-linear map from $\Hcal^*(N,N\setminus F)$ into
$\Hcal^*_F(N)$.

Let  $\beta\in \Acal^*(N\setminus F)$. If $\chi$ is a function on
$N$ which is identically $1$ on a neighborhood of $F$, note that
$d\chi\beta$ defines a differential form on $N$, since $d\chi$ is
equal to $0$ in a neighborhood of $F$.

\begin{prop}\label{alphau}
For any open neighborhood $U$ of $F$, we choose $\chi\in\f(N)$ with
support in $U$ and equal to $1$ in a neighborhood of $F$.

$\bullet$ The map
\begin{equation}\label{eq:p-U-chi}
\p_U^{\chi}\left(\alpha,\beta\right)=\chi\alpha + d\chi\beta
\end{equation}
defines a homomorphism of complexes
$\p_U^{\chi}:\Acal^*(N,N\setminus F)\to \Acal^*_U(N).$

Let $\alpha\in \Acal^*(N)$ be a \emph{closed} form
and $\beta\in \Acal^*(N\setminus F)$ such that $\alpha|_{N\setminus
F}=d\beta$. Then $\p_U^{\chi}(\alpha,\beta)$ is a \emph{closed}
differential form supported in $U$.

$\bullet$ The  cohomology class of  $\p_U^{\chi}(\alpha,\beta)$ in
$\Hcal^*_U(N)$ does not depend on $\chi$. We denote this class by
$\p_U(\alpha,\beta)\in \Hcal^*_U(N)$.

$\bullet$ For any neighborhoods $V\subset U$ of $F$, we have
$f_{U,V}\circ \p_V=\p_U$.
\end{prop}
\begin{proof}
The equation $\p_U^{\chi}\circ \,d_{\rm rel}=d \circ \p_U^{\chi}$ is
immediate to check. In particular $\p_U^{\chi}(\alpha,\beta)$ is
closed, if $d_{\rm rel}\left(\alpha,\beta\right)=0$. This
proves the first point. For two different choices $\chi$ and
$\chi'$, we have $\p_U^{\chi}(\alpha,\beta)-\p_U^{\chi'}(\alpha,\beta)
= d\left((\chi-\chi')\beta\right)$. Since $\chi-\chi'=0$ in a neighborhood of $F$,
the term  $(\chi-\chi')\beta$ is a well defined  element of
$\Acal^*_U(N)$. This proves the second  point. Finally, the last point
is immediate, since
$\p_U^{\chi}(\alpha,\beta)=\p_V^{\chi}(\alpha,\beta)$ for
$\chi\in\f(N)$ with support in $V\subset U$.
\end{proof}\bigskip

\begin{defi}\label{def-alpha-beta}
Let $\alpha\in \Acal^*(N)$ be a \emph{closed} form and $\beta\in
\Acal^*(N\setminus F)$ be such that $\alpha|_{N\setminus F}=d\beta$. We
denote by    $\p_F(\alpha,\beta)\in \Hcal^*_F(N)$ the element
defined by the sequence $\p_U(\alpha,\beta)\in \Hcal^*_U(N),\,
U\in\Fcal_F$. We have then a morphism of $\Kbb$-vector spaces
\begin{equation}\label{eq:map-p}
\p_F:\Hcal^*(N,N\setminus F)\to \Hcal^*_F(N).
\end{equation}
\end{defi}


\begin{prop}

%
%
If $F_1,F_2$ are closed subsets of $N$, then we have
\begin{equation}\label{eq:fonctoriel-p-produit}
\p_{F_1\cap F_2}(a_1\diamond a_2)= \p_{F_1}(a_1)\wedge \p_{F_2}(a_2)
\end{equation}
for any $a_k\in \Hcal^*(N,N\setminus F_k)$.


\end{prop}

\begin{proof}

Let $W$ be a neighborhood of $F_1\cap F_2$. Let $V_1,V_2$ be
respectively neighborhoods of $F_1$ and $F_2$ such that $V_1\cap
V_2\subset W$. Let $\chi_i\in\f(N)$ be supported in $V_i$ and equal
to $1$ in a neighborhood of $F_i$. Then $\chi_1\chi_2$ is
supported in $W$ and equal to $1$ in a neighborhood of $F_1 \cap
F_2$. Let $\Phi_1+ \Phi_2={\rm 1}_{N\setminus (F_1\cap F_2)}$ be a
partition of unity relative to the decomposition $N\setminus
(F_1\cap F_2)=N\setminus F_1\cup N\setminus F_2$.

Then one checks easily that
$$
\p^{\chi_1}_{V_1}(a_1)\wedge\ \p^{\chi_2}_{V_2}(a_2)-
\p^{\chi_1\chi_2}_{W}(a_1\diamond_\Phi a_2)
$$
is equal to
$$
 d\Big((-1)^{k_1+1}\chi_1 d\chi_2
(\Phi_1\beta_1\beta_2)+ (-1)^{k_1}\chi_2
d\chi_1(\beta_1\Phi_2\beta_2)\Big)
$$
for $d_{\rm rel}$-closed forms
$a_i=(\alpha_i,\beta_i)\in\Acal^{k_i}(N,N\setminus F_i)$. Remark
that $\Phi_1\beta_1\beta_2$ is defined on $N\setminus F_2$, so
that $d\chi_2 (\Phi_1\beta_1\beta_2)$ is well defined on $N$ and
supported in $V_2$. Thus the form $(-1)^{k_1+1}\chi_1 d\chi_2
(\Phi_1\beta_1\beta_2)+ (-1)^{k_1}\chi_2
d\chi_1(\beta_1\Phi_2\beta_2)$ is well defined on $N$ and
supported in $V_1\cap V_2\subset W$. This proves that
$\p_{F_1}(a_1)\wedge\p_{F_2}(a_2)=\p_{F_1\cap F_2}(a_1\diamond
a_2)$.
\end{proof}

\medskip

The map $\p_F:\Hcal^*(N,N\setminus F)\longrightarrow \Hcal^*_F(N)$ factors the
natural map \break $\Hcal^*(N,N\setminus F)\to \Hcal^*(N)$.

\subsection{Integration}\label{sec:integration}

We consider first the case where $F$ is a \emph{compact} subset of
an oriented manifold $N$. Let $\Hcal^*(N,N\setminus F)$ be the relative
$\Kbb$-cohomology group.  Let $\pi:N\to \{\bullet\}$ be the projection to the point. We will
describe an integration morphism $\pi_*: \Hcal^*(N,N\setminus F)\to \Kbb$.

We have a $\Kbb$-linear map
\begin{equation}\label{eq-map-p-compact}
    \p_c :\Hcal^*(N,N\setminus F)\longrightarrow \Hcal^*_c(N)
\end{equation}
which is equal to the composition of $\p_F$ with the natural map
$\Hcal^*_F(N)\to \Hcal^*_c(N)$. If $a\in \Hcal^*(N,N\setminus F)$
is represented by the $d_{\rm rel}$-closed differential form
$(\alpha,\beta)\in\Acal^*(N,N\setminus F)$, the class $\p_c(a)\in
\Hcal^*_c(N)$ is represented by the differential form
$\p_U^\chi(\alpha,\beta)=\chi\alpha+d\chi\beta$ where $\chi$ is a
function with \emph{compact} support.

\begin{defi}\label{int}
If $a\in \Hcal^*(N,N\setminus F)$, then $\pi_*(a)\in \Kbb$ is defined by
$$\pi_*(a):=\int_N \p_c(a).$$
\end{defi}

If $N$ is {\em compact}, the elements
$\alpha$ and $\p_c(a)$ coincide in $\Hcal^*(N)$, hence
$\pi_*(a)=\int_N \alpha$. When $N$ is {\em non-compact}, an interesting situation is the
case of a relative class $a=[\alpha,\beta]$ where the
closed form $\alpha$ is {\em integrable}. The two terms $\pi_*(a)$ and $\int_N \alpha$
are defined. However, it is usually not true
that they coincide. An interesting case is the relative
Thom form $\tur(V)$ of a real oriented vector space $V$  (see Section \ref{sec:Thom-class}). 
Here $N=V$, $F=\{0\}$, and the relative class $\tur(V)$ is represented by $[0,\beta]$ with
$\beta$ a particular closed  real form on $V\setminus \{0\}$. Here the
integral of $\alpha=0$ is equal to $0$, while $\pi_* (\tur(V))=1$. See
Example \ref{basic2}.

In some important  cases studied in  Subsection
\ref{sec:gaussian}, we will however prove  that the integral of
$\p_c(a)$ is indeed the same  as  the integral of  $\alpha$. As we have
\begin{eqnarray}\label{eq:trangression-elementaire}
\p_U^\chi(\alpha,\beta)-\alpha&=& (\chi-1)\alpha+d\chi\beta\nonumber \\
                              &=& d\Big((\chi-1)\beta\Big),
\end{eqnarray}
 the comparison between the integral of $\p_c(a)$ and the one of
$\alpha$ will follow from the careful study of the behavior on $N$
of the form $(\chi-1)\beta$.

\bigskip

We consider now the case of an oriented real vector bundle $\pi:\Vcal\to M$. We will
describe a push-forward $\Kbb$-linear map $\pi_*: \Hcal^*(\Vcal,\Vcal\setminus M)\to \Hcal^*(M)$.

Let $\Acal^*_{\cf}(\Vcal)$ be the $\Kbb$-subalgebra of $\Acal^*(\Vcal)$ formed by the
differential forms which have a compact support in the fibers of $\pi$. Let
$\Hcal^*_{\cf}(\Vcal)$ be the corresponding $\Kbb$-cohomology space. We have a morphism
$\int_{\rm fiber}:\Hcal^*_{\cf}(\Vcal)\to \Hcal^*(M)$ of integration along the fibers.

We define a $\Kbb$-linear map
\begin{equation}\label{eq-map-p-compact-fiber}
    \p_{\cf} :\Hcal^*(\Vcal,\Vcal\setminus M)\longrightarrow \Hcal^*_{\cf}(\Vcal)
\end{equation}
by setting that $\p_{\cf}([\alpha,\beta])$ is the class represented by
$\chi\alpha+d\chi\beta$, where $\chi$ is a function on $\Vcal$ with \emph{compact}
support in the fibers, and equal to $1$ in a neighborhood of the zero section.

\begin{defi}\label{int-fiber}
If $a\in \Hcal^*(\Vcal,\Vcal\setminus M)$, the class $\pi_*(a)\in \Hcal^*(M)$ is defined by
$$
\pi_*(a):=\int_{\rm fiber}\p_{\cf}(a).
$$
\end{defi}

In Section \ref{sec:Thom-class}, we will describe a relative Thom class $\tur(\Vcal)$
which is characterized by the fact that $\pi_*(\tur(\Vcal))=1$ in $\Hcal^*(M)$.

\section{Quillen's relative Chern Character}\label{sec:chrel-sigma}

In this section, we work with differential forms with {\em complex} coefficients.

\subsection{Chern form of a super-connection}

For an introduction to the Quillen's notion of super-connection,
see \cite{BGV}.

If $\Ecal$ is a  complex vector bundle on  a manifold $N$, we denote by
$\Acal^*(N,\End(\Ecal))$ the complex algebra of $\End(\Ecal)$-valued
differential forms on $N$.

Let $\nabla$ be  a connection on $\Ecal$. The curvature $\nabla^2$ of
$\nabla$ is a $\End(\Ecal)$-valued two-form on $N$.  Recall that the
Chern character of $\Ecal$ is the de Rham cohomology class of  the
closed differential form ${\rm
Chern}(\Ecal):=\tr(\exp(\frac{-\nabla^2}{2i\pi}))$. Here we simply
denote by  $\ch(\Ecal)\in \Hcal^*(N)$ the de Rham cohomology class
of $\tr(\exp(\nabla^2))$. We will call it  the (non normalized)
Chern character of $\Ecal$.

More generally, let $\Ecal=\Ecal^+\oplus \Ecal^-$ be a $\Zbb_2$-graded complex vector
bundle on a manifold $N$. Taking in account the $\Zbb_2$-grading of
$\End(\Ecal)$, the algebra $\Acal^*(N,\End(\Ecal))$ is a
$\Zbb_2$-graded algebra: for example
$[\Acal^*(N,\End(\Ecal))]^+$ is equal to \break 
$\Acal^+(N,\End(\Ecal)^+)\oplus\Acal^-(N,\End(\Ecal)^-)$.
 The super-trace  on $\End(\Ecal)$ extends to
a $\Cbb$-linear map $\str:\Acal^*(N,\End(\Ecal))\to \Acal^*(N)$.

Let $\A$ be a super-connection on $\Ecal$ and $\mathbf{F}=\A^2$ be its
curvature, an element of $[\Acal^*(N,\End(\Ecal))]^+$. {\em The Chern
form} of $(\Ecal,\A)$ is the closed differential form
$$
\ch(\A):=\str(\e^{\F}).
$$


We will use the following transgression formulaes.

\begin{prop}\label{trans}
$\bullet$ Let $\A_t$, for $t\in \Rbb$, be  a one parameter family of
super-connections on $E$, and let $\frac{d}{dt}\A_t\in
[\Acal^*(N,\End(\Ecal))]^-$. Let $\F_t$ be the  curvature of $\A_t$. Then
one has 
\begin{equation}\label{transgression}
\frac{d}{dt}\ch(\A_t)=d\left(\str\Big((\frac{d}{dt}\A_t)
\e^{\F_t}\Big)\right).
\end{equation}
$\bullet$ Let $\A(s,t)$ be a two-parameter family of
super-connections. Here $s,t\in \Rbb$. We denote by   $\F(s,t)$ the
curvature of $\A(s,t)$. Then:
\begin{eqnarray*}
\lefteqn{\frac{d}{ds}\str\Big((\frac{d}{dt}\A(s,t))\,
\e^{\F(s,t)}\Big)-
\frac{d}{dt}\str\Big((\frac{d}{ds}\A(s,t))\, \e^{\F(s,t)}\Big)} \nonumber\\
& &= d\left(\int_{0}^1\str\Big((\frac{d}{ds}\A(s,t)) \e^{u
\F(s,t)}(\frac{d}{dt}\A(s,t))\,\e^{(1-u)\F(s,t)}\Big)du\right).
\end{eqnarray*}

\end{prop}

\begin{proof} These formulae are well known, and are derived easily from the two
identities: $\F=\A^2$, and $d\str(\alpha)=\str[\A,\alpha]$ for any
$\alpha \in \Acal^*(N,\End(\Ecal))$ (see \cite{BGV}).
\end{proof}



In particular, the cohomology class defined  by $\ch(\A)$ in
$\Hcal^*(N)$ is independent on the choice of the super-connection
$\A$ on $\Ecal$. By definition, this is the Chern character $\ch(\Ecal)$
of $\Ecal$. By choosing $\A=\left(\begin{array}{cc}
\nabla^+& 0\\
0 & \nabla^-\\
\end{array}\right)$
where $\nabla^{\pm}$ are connections on $\Ecal^{\pm}$, this class is
just $\ch(\Ecal^+)-\ch(\Ecal^-)$. However, different choices of $\A$
define very different looking representatives of $\ch(\Ecal)$.

\subsection{Quillen's relative Chern character of a
morphism}\label{sec:quillen-chern-morphism}

Let $\Ecal=\Ecal^+\oplus \Ecal^-$ be a $\Zbb_2$-graded complex vector bundle on a
manifold $N$ and $\sigma: \Ecal^+ \to \Ecal^-$ be a smooth morphism.
At each point $n\in N$, $\sigma(n): \Ecal^+_n\to \Ecal^-_n$ is a
linear map. The \emph{support} of $\sigma$ is the closed subset of
$N$
$$
\supp(\sigma)=\{n\in N\mid \sigma(n)\ {\rm is\ not \ invertible}\}.
$$


Recall that the morphism $\sigma$ is {\em elliptic} when $\supp(\sigma)$ is  compact : in this situation
the data $(\Ecal^+,\Ecal^-,\sigma)$ defines an element
of the $\KK^0$-theory of $N$.

In the following, {\em we do not} assume $\sigma$ elliptic. We
recall Quillen's construction \cite{Quillen85}  of a  $\Cbb$-cohomology
class $\ch_{\rm rel}(\sigma)$ in $\Hcal^*(N,N\setminus
\supp(\sigma))$. The definition will involve several choices. We
choose Hermitian structures on $\Ecal^{\pm}$ and  a super-connection
$\A$ on $\Ecal$ {\em without $0$ exterior degree term}.

We associate to the morphism $\sigma$ the odd Hermitian endomorphism of $\Ecal$ defined by
\begin{equation}\label{eq:v-sigma}
 v_\sigma=
\left(\begin{array}{cc}
0 & \sigma^*\\
\sigma & 0\\
\end{array}\right).
\end{equation}
 Then
$ v_\sigma^2=\left(\begin{array}{cc}
\sigma^*\sigma& 0\\
0 &\sigma\sigma^*\\
\end{array}\right)
$ is a non negative even Hermitian  endomorphism of $\Ecal$. The support
of $\sigma$ coincides with the set of elements $n\in N$ where the
spectrum of  $v_\sigma^2(n)$ contains $0$.

\begin{defi}
Let $E$ be a finite dimensional Hermitian vector space.
If $H$ is an Hermitian endomorphism of $E$ and $h\in\Rbb$, we
write $H\geq h$ when $(Hw,w)\geq h\|w\|^2$ for any $w\in E$. Then
$H\geq h$ if and only if the smallest eigenvalue of $H$ is larger
than $h$.
\end{defi}

\medskip

 Consider the family of super-connections
\begin{equation}\label{eq:A-sigma}
    \A^{\sigma}(t)=\A+i t\, v_\sigma, \quad t\in \Rbb.
\end{equation}
The curvature of $\A^{\sigma}(t)$ is the even element
$\F(\sigma,\A,t)\in[\Acal^*(N,\End(\Ecal))]^+$ defined by :
\begin{equation}\label{eq:F-A-sigma}
\F(\sigma,\A,t)=(it v_\sigma+\A)^2=-t^2 v_\sigma^2+it
[\A,v_\sigma]+\A^2.
\end{equation}
Here $-t^2 v_\sigma^2$ is the term of exterior degree $0$. As the
super-connection $\A$ do not have $0$ exterior degree term, both
elements $it [\A,v_\sigma]$ and  $\A^2$ are sums of terms with
strictly positive exterior degrees. For example, if
$\A=\nabla^+\oplus \nabla^-$ is a direct sum of connections, then
$it [\A,v_\sigma]\in \Acal^1(N,\End(\Ecal)^-)$ and $\A^2\in
\Acal^2(N,\End(\Ecal)^+)$.

\begin{defi}\label{def:ch-A-sigma-t}
We denote by $\ch(\sigma,\A,t)$ the Chern form of $(\Ecal,\A^\sigma(t))$, that is
$$
\ch(\sigma,\A,t):=\str\left(\e^{\F(\sigma,\A,t)}\right).
$$
\end{defi}

As $iv_\sigma=\frac{d}{dt}\A^\sigma(t)$, we have the transgression formula
$\frac{d}{dt}\ch(\sigma,\A,t)=$ \break $-d(\eta(\sigma,\A,t))$ with
\begin{equation}\label{eq:transgression-eta}
\eta(\sigma,\A,t):=- \str\left(iv_\sigma\,\e^{\F(\sigma,\A,t)}\right).
\end{equation}
After integration, the transgression formula gives the following equality of differential forms
on $N$
\begin{equation}\label{eq:transgression-integral}
    \ch(\A)-\ch(\sigma,\A,t)=d\left(\int_0^t\eta(\sigma,\A,s)ds\right),
\end{equation}
since $\ch(\A)=\ch(\sigma,\A,0)$.

\begin{prop}\label{estimatesord}
Let $\Kcal$ be a compact subset of $N$ and let $h\geq 0$ be such that
$v_\sigma^2(n)\geq h$ when $n\in \Kcal$. There exists a polynomial
$\Pcal_{\Kcal}$ of degree $\dim N$ such that, on $\Kcal$,
\begin{equation}\label{estN}
\Big|\!\Big| \e^{\F(\sigma,\A,t)}\Big|\!\Big|\leq  \Pcal_{\Kcal}(t)
\e^{-ht^2} \quad \mathrm{for\ all\ }t\geq 0.
\end{equation}
  In particular, when $\Kcal$ is contained in
$N\setminus \supp(\sigma)$, then $\ch(\sigma,\A,t)$ and
$\eta(\sigma,\A, t)$ tends to $0$
  exponentially fast when
$t$ tends to infinity.
\end{prop}

\begin{proof}
 We work on
$\End(E)\otimes \Acal$, where $\Acal=\oplus_{k=0}^{\dim
N}\Acal^k(N)$. To estimate $\|\e^{\F(\sigma,\A,t)}\|$, we employ
Lemma \ref{suffit} of the Appendix, with $H=t^2v_\sigma^2$, and
$R=-it [\A,v_\sigma]-\A^2$. Here $R$ is a sum of $\End(E)$-valued
differential forms on $N$ with strictly positive exterior degrees.
Remark that $R$ is a polynomial in $t$ of degree $1$. Lemma
\ref{suffit} gives us the estimate $\|\e^{\F(\sigma,\A,t)}\|\leq
\Pcal(\|R\|) \e^{-ht^2}$ with $\Pcal$ an explicit polynomial with
positive coefficients of degree $\dim N$. Using the fact that
$\|R\|\leq at+b$ on $\Kcal$, we obtain the estimate (\ref{estN})
on $N$.

If $\Kcal$ is contained in $N\setminus \supp(\sigma)$, we can find
$h>0$ such that $v_\sigma^2(n)\geq h$ when $n\in \Kcal$. Thus we see
that $\|\e^{\F(\sigma,\A,t)}\|$ decreases exponentially fast, when
$t$ tends to infinity.
\end{proof}

The former estimates allows us to take the limit $t\to\infty$ in
(\ref{eq:transgression-integral}) on the open subset $N\setminus
\supp(\sigma)$. There, the differential form
$\ch(\sigma,\A,t)=\str\left(\e^{\F(\sigma,\A,t)}\right)$ tends to $0$
as $t$ goes to $\infty$, and we get the following important lemma
due to Quillen.
\begin{lemm}{\rm \cite{Quillen85}}\label{lem:quillen}
We can define on $N\setminus \supp(\sigma)$ the differential form
\begin{equation}
  \label{eq:beta}
  \beta(\sigma,\A)=\int_{0}^\infty\eta(\sigma,\A,t)dt
\end{equation}
and we have $\ch(\A)|_{N\setminus
\supp(\sigma)}=d\left(\beta(\sigma,\A)\right)$.
\end{lemm}

\medskip

We are in the situation of Definition \ref{relcoh}. The closed form
$\ch(\A)$ on $N$ and the form $\beta(\sigma,\A)$ on
$N\setminus\supp(\sigma)$ define an even relative cohomology class
$\left[\ch(\A),\beta(\sigma,\A)\right]$ in $\Hcal^*(N,N\setminus
\supp(\sigma))$.

\begin{prop}\label{inds}

$\bullet$ The class
$\left[\ch(\A),\beta(\sigma,\A)\right]\in\Hcal^*(N,N\setminus
\supp(\sigma))$ does not depend on the choice of $\A$, nor on the
Hermitian structure on $\Ecal$. We denote it by $\chr(\sigma)$ and
we call it the {\em Quillen Chern character}.

$\bullet$ Let $F$ be a closed subset of $N$. For $s\in [0,1]$, let
$\sigma_s:\Ecal^+\to \Ecal^-$ be a family of smooth morphisms such
that $\supp(\sigma_s)\subset F$. Then all classes $\chr(\sigma_s)$
coincide in $\Hcal^*(N,N\setminus F)$.
\end{prop}

\begin{proof} 
Let us prove the first point. Let $\A_s, s\in [0,1]$, be a smooth
one parameter family of super-connections on $\Ecal$ without $0$
exterior degree terms. Let $\A(s,t)=\A_s+it v_\sigma$. Thus
$\frac{d}{ds}\A(s,t)=\frac{d}{ds}\A_s$ and
$\frac{d}{dt}\A(s,t)=iv_\sigma$. Let $\F(s,t)$ be the curvature of
$\A(s,t)$. We have $\frac{d}{ds}\ch(\A_s)=d\gamma_s$ with
$\gamma_s=\str\Big((\frac{d}{ds}\A_s)\e^{\F(s,0)}\Big)$ and
$\eta(\sigma,\A_s,t)=-\str\left((\frac{d}{dt}\A(s,t))\e^{\F(s,t)}\right)$.
We apply the double transgression formula of Proposition
\ref{trans}, and we obtain 
\begin{equation}\label{z1ord}
\frac{d}{ds}\eta(\sigma,\A_s,t)=-\frac{d}{dt}
\str\Big((\frac{d}{ds}\A(s,t))\e^{\F(s,t)}\Big)-d\nu(s,t)
\end{equation}
 with
\begin{eqnarray*}
   \nu(s,t)&=& \int_{0}^1 \str\Big((\frac{d}{ds}\A(s,t)) \e^{u
\F(s,t)}(\frac{d}{dt}\A(s,t))\e^{(1-u)\F(s,t)}\Big)du\\
           &=& \int_{0}^1 \str\Big((\frac{d}{ds}\A_s)
\e^{u \F(s,t)}(iv_\sigma)\e^{(1-u)\F(s,t)}\Big)du.
\end{eqnarray*}

For $u,s\in [0,1]$ and $t\geq 0$, we consider the element of
$\Acal^*(N,\End(\Ecal))$ defined by
$I(u,s,t)= (\frac{d}{ds}\A_s) \e^{u \F(s,t)}(iv_\sigma)\e^{(1-u)\F(s,t)}$.

On a compact subset $\Kcal$  of $N\setminus F$, $\sigma$ is
invertible and we can find $h>0$ such that  $v_\sigma^2(n)\geq h$
when $n\in \Kcal$. We have $u \F(s,t)=-t^2uv_{\sigma}^2-u
R_{t,s}$, with $R_{t,s}=-it [\A_s,v_\sigma]-\A_s^2$ which is a sum
of terms with strictly positive exterior degrees. Remark that
$R_{t,s}$ is a polynomial of degree $1$ in $t$. By the estimate of
Lemma \ref{suffit} of the Appendix, we obtain
\begin{eqnarray*}
  \Big|\!\Big|I(u,s,t)\Big|\!\Big| & \leq & \Big|\!\Big|\frac{d}{ds}\A_s\Big|\!\Big|\,
  \Big|\!\Big| v_\sigma \Big|\!\Big|\,
  \Big|\!\Big| \e^{u \F(s,t)}\Big|\!\Big|\,
  \Big|\!\Big| \e^{(1-u) \F(s,t)}\Big|\!\Big|\\
  &\leq & \Big|\!\Big|\frac{d}{ds}\A_s\Big|\!\Big|\,
  \Big|\!\Big| v_\sigma \Big|\!\Big|\,
  \e^{-ht^2}\Pcal(u\|R_{t,s}\|)
  \Pcal((1-u)\|R_{t,s}\|)
\end{eqnarray*}
where $\Pcal$ is a polynomial of degree less or equal to $\dim N$.
So, we can find a constant $C$ such that $\|I(u,s,t)\|\leq
C(1+t^2)^{\dim N}\e^{-ht^2}$ for all $u,s\in [0,1]$ and $t\geq 0$.
Thus we can integrate Equation (\ref{z1ord}) in $t$, from 0 to
$\infty$. Since $-\int_{0}^{\infty}\frac{d}{dt}[
\str((\frac{d}{ds}\A(s,t))\e^{\F(s,t)})]dt$ $=\gamma_s$, it follows
that
\begin{equation}\label{second}
\frac{d}{ds}\ch(\A_s)=d\gamma_s,\hspace{1cm}\frac{d}{ds}\beta(\sigma,\A_s)
=\gamma_s-d\epsilon_s
\end{equation}
where
$\epsilon_s=\int_{0}^{\infty}\nu(s,t)dt$ and $\nu(s,t)= \int_{0}^{1}I(u,s,t)du$.
The first equality in Equations (\ref{second}) holds on $N$, and the
second  on  $N\setminus \supp(\sigma)$. These  equations
(\ref{second}) exactly mean that
$$
\frac{d}{ds}\left(\ch(\A_s),\beta(\sigma,\A_s)\right)=d_{\rm
rel}\left(\gamma_s,\epsilon_s\right).
$$
So the cohomology class
$\left[\ch(\A),\beta(\sigma,\A)\right]$ in $\Hcal^*(N,N\setminus
\supp(\sigma))$
 does not depend on the choice of $\A$. With a similar proof, we see that this cohomology class
does not depend  on the choice of Hermitian structure on $\Ecal$.

The proof of the second point is similar.
\end{proof}

We have defined a representative of the relative Chern class
$\chr(\sigma)$ using the one-parameter family $\A^\sigma(t)$ of
super-connections, for $t$ varying between $0$ and $\infty$. We
can define another representative as follows. We have
$\ch(\sigma,\A,t)=d(\beta(\sigma,\A,t))$ with
\begin{equation}\label{eq:beta-sigma-1}
    \beta(\sigma,\A,t)=\int_t^\infty\eta(\sigma,\A,s)ds.
\end{equation}
\begin{lemm}\label{lem:other-representant}
 For any $t\in\Rbb$, the relative Chern character $\chr(\sigma)$ satisfies
$\chr(\sigma)=[\ch(\sigma,\A,t),\beta(\sigma,\A,t)]$ in $\Hcal^*(N,N\setminus\supp(\sigma))$.
\end{lemm}

\begin{proof} It is easy to check that
\begin{equation}\label{eq-chrel-t-delta}
\Big(\ch(\A),\beta(\sigma,\A)\Big)-\Big(\ch(\sigma,\A,t),\beta(\sigma,\A,t)\Big)
=d_{\rm rel}(\delta(t),0),
\end{equation}
with $\delta(t)=\int_0^t\eta(\sigma,\A,s)ds$.
\end{proof}

\begin{rema}
 Quillen relative Chern character seems to be very related to the
``multiplicative $K$-theory" defined by Connes-Karoubi (see
\cite{kar1}, \cite{kar2}). For example, even if $\Ecal^+, \Ecal^-$
are flat bundles, the Quillen Chern character is usually non zero,
as it also encodes  the odd closed differential form
$\omega=\beta(\sigma,\A)$.
\end{rema}

\section{Multiplicative property of $\chr$}\label{sec:chr-multiplicative}

Let $\Ecal_1, \Ecal_2$ be two  $\Zbb_2$-graded complex vector bundles on a
manifold $N$. The space $\Ecal_1\otimes \Ecal_2$ is a
$\Zbb_2$-graded complex vector bundle with even part $\Ecal_1^+\otimes
\Ecal_2^+\oplus \Ecal_1^-\otimes \Ecal_2^-$ and odd part
$\Ecal_1^-\otimes \Ecal_2^+\oplus \Ecal_1^+\otimes \Ecal_2^-$.


The complex super-algebra $\Acal^*(N,\End(\Ecal_1\otimes \Ecal_2))$ can be
identified with \break
$\Acal^*(N,\End(\Ecal_1))\otimes\Acal^*(N,\End(\Ecal_2))$ where
the tensor is taken in the sense of super-algebras. Then, if
$v_k\in \Acal^0(N,\End(\Ecal_k)^-)$ are odd endomorphisms, we have
$(v_1\otimes{\rm Id}_{\Ecal_2}+ {\rm Id}_{\Ecal_1}\otimes
v_2)^2=(v_1)^2\otimes {\rm Id}_{\Ecal_2}+ {\rm
Id}_{\Ecal_1}\otimes (v_2)^2$.

\bigskip

Let $\sigma_1:\Ecal_1^+\to \Ecal_1^-$ and $\sigma_2:\Ecal_2^+\to
\Ecal_2^-$ be two smooth morphisms. With the help of Hermitian
structures, we define the morphism
$$
\sigma_1\odot \sigma_2: \left(\Ecal_1\otimes \Ecal_2\right)^+
\longrightarrow \left(\Ecal_1\otimes \Ecal_2\right)^-
$$
by $\sigma_1\odot \sigma_2:= \sigma_1\otimes {\rm
Id}_{\Ecal_2^+}+{\rm Id}_{\Ecal_1^+}\otimes \sigma_2+{\rm
Id}_{\Ecal_1^-}\otimes\sigma_2^*+ \sigma_1^*\otimes {\rm
Id}_{\Ecal_2^-}$.

\bigskip

Let $v_{\sigma_1}$ and $v_{\sigma_2}$ be the odd Hermitian endomorphisms of $\Ecal_1,\Ecal_2$
associated to $\sigma_1$ and $\sigma_2$ (see (\ref{eq:v-sigma})). Then $v_{\sigma_1\odot
\sigma_2}= v_{\sigma_1}\otimes {\rm Id}_{\Ecal_2}+{\rm
Id}_{\Ecal_1}\otimes v_{\sigma_2}$ and $v_{\sigma_1\odot
\sigma_2}^2=v_{\sigma_1}^2\otimes {\rm Id}_{\Ecal_2}+{\rm
Id}_{\Ecal_1}\otimes v_{\sigma_2}^2$. It follows
that
$$
\supp(\sigma_1\odot \sigma_2)= \supp(\sigma_1)\cap \supp(\sigma_2).
$$



We can now state one of the main result of this paper.

\begin{theo}\label{theo:chrel-produit}{\rm \bf (Quillen's Chern character is
multiplicative)} Let $\sigma_1,\sigma_2$ be two morphisms over $N$.
The relative cohomology classes
\begin{itemize}
  \item $\chr(\sigma_k)\in\Hcal^*(N,N\setminus\supp(\sigma_k))$,
  $k=1,2$,
  \item
  $\chr(\sigma_1\odot\sigma_2)\in
  \Hcal^*(N,N\setminus(\supp(\sigma_1)\cap\supp(\sigma_2)))$,
\end{itemize}
satisfy the following equality
$$
\chr(\sigma_1\odot\sigma_2)=\chr(\sigma_1)\diamond\chr(\sigma_2)
$$
in $\Hcal^*(N,N\setminus(\supp(\sigma_1)\cap\supp(\sigma_2)))$.
Here $\diamond$ is the product of relative classes (see
(\ref{eq:produit-relatif})).
\end{theo}

\begin{proof}
For $k=1,2$, we choose  super-connections $\A_k$, without $0$
exterior degree terms. We consider the closed forms
$c_k(t):=\ch(\sigma_k,\A_k,t)$ and the transgression forms
$\eta_k(t):=\eta(\sigma_k,\A_k,t)$ so that
$\frac{d}{dt}(c_k(t))=-d(\eta_k(t))$. Let
$\beta_k=\int_0^{\infty}\eta_k(t)dt.$ A representative of
$\chr(\sigma_k)$ is $(c_k(0),\beta_k)$.

For the symbol $\sigma_1\odot\sigma_2$, we consider $\A(t)=\A+it
v_{\sigma_1\odot \sigma_{2}}$ where $\A=\A_1\otimes
\id_{\Ecal_2}+\id_{\Ecal_1}\otimes \A_2$. Then
$\ch(\A)=c_1(0)c_2(0)$. Furthermore, it is easy to see that the
transgression form  for the family $\A(t)$ is
$\eta(t)=\eta_1(t)c_2(t)+c_1(t)\eta_2(t).$ Let
$\beta_{12}=\int_{0}^{\infty}\eta(t)dt.$ A representative of
$\chr(\sigma_1\odot\sigma_2)$ is $(c_1(0)c_2(0),\beta_{12})$.

We consider the open subsets
$U=N\setminus(\supp(\sigma_1)\cap\supp(\sigma_2))$, and
$U_k=N\setminus\supp(\sigma_k)$. Let $\Phi_1+\Phi_2= {\rm 1}_{U}$
be a partition of unity subordinate to the decomposition
$U=U_1\cup U_2$. The proof will be completed if one shows that
\begin{equation}\label{prodrel}
(c_1(0),\beta_1)\diamond_\Phi(c_2(0),\beta_2)-(c_1(0)c_2(0),\beta_{12})
\end{equation}
is $d_{\rm rel}$-exact. We need the following lemma.
\begin{lemm}\label{prop:B-defined}
 The integrals
\begin{eqnarray*}
  B_1 &=& \int\!\!\!\int_{0\leq t \leq  s}\Phi_1  \eta_1(s) \wedge \eta_2(t)ds \,dt, \\
  B_2 &=&  \int\!\!\!\int_{0\leq s \leq t} \Phi_2 \eta_1(s) \wedge \eta_2(t)ds \,dt
\end{eqnarray*}
are well defined differential forms on $U$.
\end{lemm}

\begin{proof}
The function $(s,t)\mapsto \Phi_1  \eta_1(s) \wedge \eta_2(t)$ is
a function on $\Rbb^2$ with values in $\Acal^*(U)$. We have to see
that the integral $B_1$ is convergent on the domain  $0\leq t \leq
s$. This fact follows directly from the estimates of Proposition
  \ref{estimatesord}. Indeed, let $\Kcal$ be a compact subset of $U$. Since
$\Phi_1$ is supported on $U_1=N\setminus \supp(\sigma_1)$, there
exists $h>0$, and a polynomial $\Pcal_1$ in $s$ such that, on
$\Kcal$,
$$
\|\Phi_1 \eta_1(s)\| \leq  \Pcal_1(s) \e^{-hs^2}\quad \mathrm{for}
\quad s\geq 0.
$$
On the other hand, there exists a polynomial $\Pcal_2$ in $t$ such
that, on $\Kcal$,
$$
\| \eta_2(t)\| \leq  \Pcal_2(t)\quad \mathrm{for} \quad t\geq 0.
$$
Then, when $0\leq t \leq  s$, we have, on $\Kcal$: $ \|\Phi_1
\eta_1(s) \wedge \eta_2(t)\|\leq \Pcal_1(s)\Pcal_2(t)\, \e^{-hs^2}
$ and the integral $B_1$ is absolutely convergent on $0\leq t< s$.
Reversing the role $1\leftrightarrow 2$, we prove in the same way
that $B_2$ is well defined.
\end{proof}

We now prove that (\ref{prodrel}) is equal to $d_{\rm rel}\Big(0,B_1-B_2\Big)$.
Indeed
$$
(c_1(0),\beta_1)\diamond_\Phi(c_2(0),\beta_2)=\Big(c_1(0)c_2(0),\Phi_1 \beta_1
c_2(0)+c_1(0)\Phi_2\beta_2-d\Phi_1 \beta_1\beta_2\Big),
$$
so that (\ref{prodrel}) is equal to $\Big(0,\Phi_1 \beta_1 c_2(0)+c_1(0)\Phi_2\beta_2-d\Phi_1
\beta_1\beta_2-\beta_{12}\Big)$.
Thus we need to check that $ dB_2-dB_1=\Phi_1 \beta_1
c_2(0)+c_1(0)\Phi_2\beta_2-d\Phi_1 \beta_1\beta_2-\beta_{12}$. We
have $dB_1=R_1+S_1$ with
\begin{eqnarray*}
  R_1 &=& d\Phi_1 \wedge \int\!\!\!\int_{0\leq t \leq s}\eta_1(s)\eta_2(t)\, ds\, dt, \\
  S_1 &=& \Phi_1\int\!\!\!\int_{0\leq t \leq  s}
 \Big(d\eta_1(s)\eta_2(t)-\eta_1(s)d\eta_2(t)\Big) ds \,dt.
\end{eqnarray*}
Now we use $d\eta_j(s)=-\frac{d}{ds}c_j(s)$, so that we obtain
\begin{eqnarray*}
  S_1 &=& \Phi_1\int\!\!\!\int_{0\leq t \leq  s}
\left((-\frac{d}{ds}c_1(s))\eta_2(t)+\eta_1(s)(\frac{d}{dt} c_2(t))\right) ds
\,dt \\
      &=& \Phi_1
\Big(\int_{0}^{\infty}
c_1(t)\eta_2(t)dt+\int_0^{\infty}\eta_1(s)c_2(s)ds\Big)-
c_2(0)\Phi_1\beta_1,
\end{eqnarray*}
that is $S_1=\Phi_1 \beta_{12}- c_2(0)\Phi_1\beta_1$. Similarly, we
compute that $dB_2$ is equal to $d\Phi_2 \wedge
\int\!\!\!\int_{0\leq s \leq t}\eta_1(s) \eta_2(t)\, ds\, dt -\Phi_2\beta_{12} +  c_1(0)\Phi_2\beta_2$.
So finally, as $d\Phi_1=-d\Phi_2$, we get
$$
dB_2-dB_1=-d\Phi_1
\int_{0}^\infty\!\!\!\int_{0}^\infty\eta_1(s)\eta_2(t)ds
dt-\beta_{12}+c_2(0)\Phi_1 \beta_1 + c_1(0)\Phi_2\beta_2
$$
which was the equation to prove.
\end{proof}

\section{Chern character of a morphism}\label{sec:chern-character}

We employ notations of Section \ref{sec:quillen-chern-morphism}. We work here with
differential forms with {\em complex} coefficients.

\subsection{The Chern Character}

Let $\sigma:\Ecal^+\to\Ecal^-$ be a morphism on $N$.
Following  Subsection \ref{sec:cohomologie-support}, we consider the
image of $\chr(\sigma)$ through the map $\Hcal^*(N,N\setminus
\supp(\sigma))\to \Hcal^*_{\supp(\sigma)}(N)$. Applying Propositions
\ref{alphau} and \ref{inds}, we obtain the following theorem.

\begin{theo}\label{theo:chgood}
$\bullet$ For any neighborhood $U$ of $\supp(\sigma)$, take
$\chi\in\f(N)$ which is equal to 1 in a neighborhood of
$\supp(\sigma)$ and with support contained in $U$. The differential
form
\begin{equation}\label{sigmaA}
c(\sigma,\A,\chi)=\chi\, \ch(\A) + d\chi \,\beta(\sigma,\A)
\end{equation}
is {\em closed} and supported in $U$.
 Its cohomology class $c_U(\sigma)\in \Hcal^*_U(N)$ does not depend on
the choice of $\A,\chi$ and the Hermitian structures on
$\Ecal^{\pm}$. Furthermore, the inverse family $c_U(\sigma)$ when
$U$ runs over the neighborhoods of $\supp(\sigma)$ defines a class
$$
\chg(\sigma)\in \Hcal^*_{\supp(\sigma)}(N).
$$
The image of this class in $\Hcal^*(N)$ is the  Chern character
$\ch(\Ecal)$ of $\Ecal$.

$\bullet$ Let $F$ be a closed subset of $N$. For $s\in [0,1]$, let
$\sigma_s:\Ecal^+\to \Ecal^-$ be a family of smooth morphisms such
that $\supp(\sigma_s)\subset F$. Then all classes $\chg(\sigma_s)$
coincide in $\Hcal^*_F(N)$.
\end{theo}

Using Lemma \ref{lem:other-representant} we get
\begin{lemm}\label{lem:c-sigma-t}
For any $t\geq 0$, the class $\chg(\sigma)$ can be defined with the
forms $c(\sigma,\A,\chi,t)=\chi\, \ch(\sigma,\A,t) +
d\chi\,\beta(\sigma,\A,t)$.
\end{lemm}
\begin{proof} It is due to the following transgression
\begin{equation}\label{eq:trangression-c-t}
    c(\sigma,\A,\chi)-c(\sigma,\A,\chi,t)=d(\chi\delta(t)),
\end{equation}
which follows from (\ref{eq-chrel-t-delta}).
\end{proof}

In some situations, {\em Quillen's Chern character} $\chq(\sigma)=\ch(\sigma,\A,1)$ enjoys good
properties relative to the integration. So it is natural to compare
the differential forms  $c(\sigma,\A,\chi)$ and $\chq(\sigma)$.

\begin{lemm}\label{retardord}
We have
$$c(\sigma,\A,\chi)-\chq(\sigma)=
d\left(\chi\int_{0}^1\eta(\sigma,\A,s)ds\right)+ d\Big((\chi-1)
\beta(\sigma,\A,1)\Big).$$
\end{lemm}
\begin{proof} This follows immediately from the transgressions
(\ref{eq:trangression-elementaire}) and (\ref{eq:trangression-c-t}).
\end{proof}

\begin{defi}\label{def:ch-compact}
    When $\sigma$ is elliptic, we denote by
\begin{equation}\label{eq:ch-c}
    \ch_c(\sigma)\in \Hcal^*_c(N)
\end{equation}
the cohomology class with compact support which is the image of
$\chr(\sigma)$ through the
map $\p_c:\Hcal^*(N,N\setminus \supp(\sigma))\to \Hcal^*_{c}(N)$ (see (\ref{eq-map-p-compact})).
\end{defi}

A representative of $\ch_c(\sigma)$ is given by $c(\sigma,\A,\chi)$,
where $\chi\in\f(N)$ is chosen with a compact support, and equal to
$1$ in a neighborhood of $\supp(\sigma)$ and $c(\sigma,\A,\chi)$ is
given by Formula (\ref{sigmaA}).

\bigskip

We will now rewrite Theorem \ref{theo:chrel-produit} for the Chern classes
$\chg$ and $\ch_c$. Let $\sigma_1:\Ecal_1^+\to \Ecal_1^-$ and
$\sigma_2:\Ecal_2^+\to \Ecal_2^-$ be two
smooth morphisms. Let $\sigma_1\odot \sigma_2: \left(\Ecal_1\otimes \Ecal_2\right)^+
\to \left(\Ecal_1\otimes \Ecal_2\right)^-$ be  their product.

Following (\ref{eq:produit-support}), the product of the elements
$\chg(\sigma_k)\in \Hcal^*_{\supp(\sigma_k)}(N)$ for $k=1,2$ belongs
to $\Hcal^*_{\supp(\sigma_1)\cap\supp(\sigma_2)}(N)=
\Hcal^*_{\supp(\sigma_1\odot\sigma_2)}(N)$.

\begin{theo}\label{theo:produit-chg}
$\bullet$ We have the equality
$$
\chg(\sigma_1)\wedge \chg(\sigma_2)=
\chg(\sigma_{1}\odot\sigma_2)\quad {\rm in}\quad  \Hcal^*_{\supp(\sigma_1\odot
\sigma_2)}(N).
$$

$\bullet$ If the morphisms $\sigma_1,\sigma_2$ are elliptic, we have
$$
\ch_c(\sigma_1)\wedge \ch_c(\sigma_2)=
\ch_c(\sigma_{1}\odot\sigma_2)\quad {\rm in} \quad \Hcal^*_{c}(N).
$$
\end{theo}

\begin{proof} The second point is a consequence of the first point.
Theorem \ref{theo:chrel-produit} tells us that
$\chr(\sigma_1\odot\sigma_2)=\chr(\sigma_1)\diamond\chr(\sigma_2)$
holds in
$\Hcal^*(N,N\setminus(\supp(\sigma_1)\cap\supp(\sigma_2)))$. We
apply now the morphism ``$\p$" and use the relation
(\ref{eq:fonctoriel-p-produit}) to get the first point.
\end{proof}

The second point of Theorem \ref{theo:produit-chg} has the
following interesting refinement. Let $\sigma_1,\sigma_2$ be two
morphisms on $N$ which are {\bf not elliptic}, and assume that the
product $\sigma_1\odot\sigma_2$ is {\bf elliptic}. Since
$\supp(\sigma_1)\cap\supp(\sigma_2)$ is compact, we consider
neighborhoods $U_k$ of $\supp(\sigma_k)$ such that
$\overline{U_1}\cap\overline{U_2}$ is compact. Let
$\chi_k\in\f(N)$ be  supported on $U_k$ and equal to $1$ in a
neighborhood of $\supp(\sigma_k)$.  Then, the  differential form
$c(\sigma_1,\A_1,\chi_1)\wedge c(\sigma_2,\A_2,\chi_2)$ is
compactly supported on $N$, and we have
$$
\ch_c(\sigma_1\odot\sigma_2)=\Big[c(\sigma_1,\A_1,\chi_1)\wedge
c(\sigma_2,\A_2,\chi_2)\Big]\quad
{\rm in} \quad \Hcal_c^*(N).
$$
Note that the differential forms $c(\sigma_k,\A_k,\chi_k)$ are not compactly supported.

\subsection{Comparison with other constructions.}

\subsubsection{Trivialization outside $\supp(\sigma)$}

Outside the support of $\sigma$, the complex vector bundles $\Ecal^+$ and
$\Ecal^-$ are ``the same", so that it is natural to construct
representatives of $\ch(\Ecal)=\ch(\Ecal^+)-\ch(\Ecal^-)$ which are
zero ``outside'' the support of $\sigma$ by the following
identifications of bundles with connections. For simplicity, we
assume in this section that $\sigma$ is elliptic.

A connection $\nabla=\nabla^+\oplus \nabla^-$ is said ``adapted'' to
the morphism $\sigma$ when the following holds
\begin{eqnarray}\label{eq:adapted}
    \nabla^- \circ \sigma    +   \sigma \circ \nabla^+ &=& 0,\\
\nabla^+ \circ \sigma^*  + \sigma^* \circ \nabla^-     &=&
0\nonumber,
\end{eqnarray}
outside a compact neighborhood of $\supp(\sigma)$. An adapted
connection is  denoted by $\nabla^{\textrm{\tiny adap}}$. It  is
easy to construct an adapted connection.


\begin{prop}\label{prop:ch=ch}
    Let $\nabla^{\textrm{\tiny adap}}$ be a connection {\em adapted} to
    $\sigma:\Ecal^+\to \Ecal^-$. Then the differential form
    $\ch(\nabla^{\textrm{\tiny adap}})$ is
    compactly supported and its cohomology class coincides with $\ch_{c}(\sigma)$
    in $\Hcal^*_c(N)$.
\end{prop}
\begin{proof} Suppose that $\nabla^{\textrm{\tiny adap}}$ satisfies
(\ref{eq:adapted}) outside a compact neighborhood $C$ of
$\supp(\sigma)$. We verify that the forms $\ch(\nabla^{\textrm{\tiny
adap}})$ as well as the form $\beta(\sigma,\nabla^{\textrm{\tiny
adap}})$ are supported on $C$. Thus if $\chi\in\f(N)$ is equal to
$1$ on $C$, we see that the differential forms
$c(\sigma,\nabla^{\textrm{\tiny adap}},\chi)$ and
$\ch(\nabla^{\textrm{\tiny adap}})$ coincide.
\end{proof}

\begin{rema}\label{Atiyahetal} If $F$ is {\em compact}, the closed differential form
$\ch(\nabla^{\textrm{\tiny adap}})$  represents the Chern
character of a difference bundle $[\Dcal^+]-[\Dcal^-]$, where
$[\Dcal^+]$ and $[\Dcal^-]$ are complex vector bundles (isomorphic outside
$F$) on a compactification $N_F$ of $N$ (see for example
\cite{B-V.92}). Thus $ \ch(\nabla^{\textrm{\tiny adap}})$ is a
representative of the Chern character as defined by Atiyah and al.
in \cite{Atiyah-Hirz61,Atiyah-Singer-2}. In this case, Theorem
\ref{theo:produit-chg} is just the multiplicativity property of
the Chern character  in  absolute theories.
\end{rema}

\subsubsection{Gaussian look}\label{sec:gaussian}

In  \cite{Mathai-Quillen}, Mathai-Quillen  gives an explicit
representative with ``Gaussian look''  of the Bott class of  a
complex vector bundle $N\to B$. The purpose of this paragraph is to
compare the Mathai-Quillen construction  of Chern characters with
``Gaussian look'' and the Quillen relative construction.

Let $N$ be a real vector bundle over a manifold $B$. We denote by
$\pi: N\to B$ the projection. We denote by $(x,\xi)$ a point of $N$
with $x\in B$ and $\xi\in N_x$.
 Let $\Ecal^\pm\to B$ be two Hermitian vector bundles. We
consider a morphism $\sigma: \pi^*\Ecal^+\to \pi^*\Ecal^-$.

We choose a metric on the fibers of the fibration $N\to B$. We work
under the following assumption on $\sigma$.

\begin{assu}\label{assum:dec-rap-ord}
The morphism $\sigma: \pi^*\Ecal^+\to \pi^*\Ecal^-$ and all its
partial derivatives have at most a polynomial growth along the
fibers of $N\to B$. Moreover we assume  that, for any compact subset
$\Kcal_B$ of $B$, there exist $R\geq 0$ and $c>0$ such
that\footnote{This inequality means that $\|\sigma(x,\xi)w\|^2\geq
c\|\xi\|^2\|w\|^2$ for any $w\in \Ecal_x$.} $v_\sigma^2(x,\xi)\geq
c\|\xi\|^2$ when $\|\xi\|\geq R$ and $x\in \Kcal_B$.
\end{assu}

We may define the sub-algebra $\Acal^*_{\dr}(N)$ of forms on $N$
such that all partial derivatives are rapidly decreasing along the
fibers. Let $\Hcal^*_{\dr}(N)$ be the corresponding cohomology
algebra. Under Assumption \ref{assum:dec-rap-ord}, the support of
$\sigma$ intersects the fibers of $\pi$  in compact sets. We have then a
canonical map from $\Hcal^*_{\supp(\sigma)}(N)$ into
$\Hcal^*_{\dr}(N)$. We will now compute the image of $\chg(\sigma)$
under this map.

Let $\nabla=\nabla^+\oplus \nabla^-$ be a connection on $\Ecal\to
B$, and consider the super-connection  $\A=\pi^*\nabla$ so that
$\A^\sigma(t)=\pi^*\nabla+i tv_\sigma$. Then, the {\em Quillen Chern character} form
$\chq(\sigma):=\ch(\sigma,\A,1)$ has a ``Gaussian'' look.
\begin{lemm}\label{prop:ch=dec-rapid}
    The differential forms $\ch(\sigma,\A,1)$ and $\beta(\sigma,\A,1)$ are
    rapidly decreasing along the fibers.
\end{lemm}

\begin{proof} The curvature of $\A^\sigma(t)$ is
$$
\F(t)= \pi^*\F -t^2v_\sigma^2+it[\pi^*\nabla,v_\sigma].
$$
Here $\F\in\Acal^2(B,\End(\Ecal))$ is the curvature of $\nabla$.
Assumption \ref{assum:dec-rap-ord} implies that
$[\pi^*\nabla,v_\sigma]\in\Acal^1(N,\End(\pi^*\Ecal))$ has at most a
polynomial growth along the fibers. Furthermore, for any compact
subset $\Kcal_B$ of the basis, there exists $R\geq 0$ and $c>0$ such
that $v_\sigma^2(x,\xi)\geq c\|\xi\|^2$ when $\|\xi\|\geq R$ and
$x\in \Kcal_B$.

\medskip

To estimate $\e^{\F(t)}$, we apply Lemma \ref{suffit} of the
Appendix, with $H=t^2 v_{\sigma}^2,$ and $R=-\pi^*\F-it
[\pi^*\nabla, v_\sigma]$. The smallest eigenvalue of $H$ is
greater or equal to
 $t^2c\|\xi\|^2$, when $\|\xi\|\geq R$, and  $R$ is a sum of terms with  strictly
positive exterior degrees. Remark that $R$ is a polynomial in $t$
of degree $1$ and is bounded in norm by a polynomial in $\|\xi\|$
along the fibers. It follows that from Lemma \ref{suffit}  that,
for $t\geq 0$, we have
$$
\Big|\!\Big| \e^{\F(t)}\Big|\!\Big|(x,\xi) \leq \Pcal(\|R\|)
\e^{-t^2c\|\xi\|^2}.
$$

Our estimates on the polynomial growth of $R$ in $t$ and $\|\xi\|$
implies that there exists a polynomial $\Qcal$ such that, for
$t\geq 0$,
\begin{equation}\label{eq:maj-F-t}
    \Big|\!\Big| \e^{\F(t)}\Big|\!\Big|(x,\xi) \leq \Qcal(t\|\xi\|)
    \e^{-t^2c\|\xi\|^2},
\end{equation}
for $(x,\xi)\in N$, $x\in \Kcal_B$, $\|\xi\|\geq R$.

This implies that $\ch(\sigma,\A,1)=\str(\e^{\F(1)})$ is rapidly
decreasing along the fibers. Consider now $\beta(\sigma,\A,1)=
-i\int_{1}^\infty \str(v_\sigma \e^{\F(t)})dt$ which is defined (at
least) for $\|\xi\|\geq R$. The estimate (\ref{eq:maj-F-t}) shows
also that $\beta(\sigma,\A,1)$ is rapidly decreasing along the
fibers. We can prove in the same way that all partial derivatives of
$\ch(\sigma,\A,1)$ and $\beta(\sigma,\A,1)$ are rapidly decreasing
along the fibers: hence $\ch(\sigma,\A,1)\in\Acal^*_{\dr}(N)$ and
$\beta(\sigma,\A,1)\in \Acal^*_{\dr}(N\setminus\supp(\sigma))$.
\end{proof}

\begin{prop}\label{prop:chg=ch-dec-rapid}
Quillen's Chern character form $\chq(\sigma)\in \Acal^*_{\dr}(N)$ represents
the image of the class $\chg(\sigma)\in\Hcal^*_{\supp(\sigma)}(N)$ in
$\Hcal^*_{\dr}(N)$.
\end{prop}
\begin{proof}
Choosing $\chi$ supported on $\|\xi\|\leq R+1$ and equal to $1$ in
a neighborhood of $\|\xi\|\leq R$, the transgression formula of
Lemma \ref{retardord}: $
c(\sigma,\A,\chi)-\ch(\sigma,\A,1)=d(\chi\int_{0}^1\eta(\sigma,\A,s)ds)+
d\left((\chi-1) \beta(\sigma,\A,1)\right)$ implies our
proposition, since the form $c(\sigma,\A,\chi)$ represents
$\chg(\sigma)$ in $\Hcal^*_{\dr}(N)$.
\end{proof}

When the fibers of $\pi: N\to B$ are oriented, we have an
integration morphism $\int_{\rm fiber}: \Hcal^*_{\dr}(N)\to \Hcal^*(B)$.

\begin{coro}
We have $\int_{\rm fiber}\chq(\sigma)=\int_{\rm fiber}\chg(\sigma)$ in $\Hcal^*(B)$.
\end{coro}


\subsection{Examples}

If $\Ecal$ is a trivial bundle $N\times E$ on a manifold $N$, an
endomorphism of $\End(\Ecal)$ is determined by a map from $N$ to
$\End(E)$. We employ the same notation for both objects, so that if
$\sigma$ is a map from $N$ to $\End(E)$, we also denote by $\sigma$
the bundle map $\sigma[n,v]=[n,\sigma(n)v]$, for $n\in N$ and $v\in
E$.

We will use  the following convention. Let $E=E^+\oplus E^-$ be a
$\Zbb_2$-graded finite dimensional complex vector space. Let $\Acal$ be a super-commutative algebra (the ring
of differential forms on a manifold for example). The elements of
the super-algebra $\Acal\otimes \End(E)$ will be  represented by
 matrices with coefficients in $\Acal$. This algebra
operates on the space $\Acal\otimes E$.  We take the following
convention: the forms are always considered as operating first:
for example, if $E^+=\Cbb$ and $E^-=\Cbb$, the matrix
$\left(\begin{array}{cc} 0 & \alpha\\\beta & 0
\end{array}\right)$
represents the operator
\begin{equation}\label{eq:convention}
\left(\begin{array}{cc} 0 & \alpha\\
\beta & 0
\end{array}\right):=\left(\begin{array}{cc} 0 & 1\\0 & 0
\end{array}\right) \alpha+
\left(\begin{array}{cc} 0 & 0\\1 & 0
\end{array}\right) \beta
\end{equation}
on $\Acal\otimes E$.

\subsubsection{The cotangent bundle $\T^*S^1$.}\label{basic1}

 We consider $\T^*S^1:=S^1\times \Rbb$ the cotangent
bundle to the circle $S^1$. The group $\KK^0(\T^*S^1)$ of K-theory
is generated, as a $\Zbb$-module, by the class $[\sigma]$ of the
following elliptic symbol.

Take $\Ecal^+=\Ecal^-$ the trivial bundles $\T^*S^1\times \Cbb$
over $\T^*S^1$. Let $u\in\f(\Rbb)$ be a function satisfying
$u(\xi)=1$ if $|\xi|>1$ and $u(\xi)=0$ if $|\xi|<1/2$. The symbol
$\sigma: \Ecal^+\to \Ecal^-$ is defined by the map $\sigma:
\T^*S^1 \to \End_\Cbb(\Cbb)=\Cbb$ :
$$
 \sigma(\e^{i\theta},\xi)=\begin{cases}
   u(\xi) \e^{i\theta}, & \text{if}\quad \xi\geq 0;\\
    u(\xi), & \text{if}\quad
   \xi\leq 0.
\end{cases}
$$
Here $\supp(\sigma)=\{(\e^{i\theta},\xi);\,u(\xi)=0\}$ is compact. Note that the
class $[\sigma]\in \KK^0(\T^*S^1)$ does not depend on the choice
of the function $u$.

We choose on $\Ecal^\pm$ the trivial connections
$\nabla^+=\nabla^-=d$ and we let $\A=\nabla^+\oplus \nabla^-$ be the trivial connection on
$\Ecal^{+}\oplus \Ecal^{-}$. Then $\ch(\A)=0$.
The curvature $\F(\sigma,\A,t)$ of the
super-connection
$$\A^\sigma(t)= \left(\begin{array}{cc} d & 0\\ 0 & d\\
\end{array}\right)+ \left(\begin{array}{cc} 0 & it u( \xi) \e^{-i\theta}\\
itu(\xi) \e^{i\theta} & 0\\ \end{array}\right)
$$
is represented by
the matrix $\left(\begin{array}{cc} a & -\overline{b}\\ b & a\\
\end{array}\right)$ where $a(t,(\e^{i\theta},\xi))=-t^2u(\xi)^2$ and
$$
b(t,(\e^{i\theta},\xi))=\begin{cases}
   -it \e^{-i\theta}(u'(\xi)d\xi-iu(\xi)d\theta), & \text{if}\quad \xi\geq 0;\\
    -it u'(\xi)d\xi, & \text{if}\quad
   \xi\leq 0.
\end{cases}
$$

Then $\e^{\F(\sigma,\A,t)}$ is represented by the
matrix $\e^{-t^2u(\xi)^2}
\left(\begin{array}{cc}
A & -\overline{b}\\b& \overline{A}\\
\end{array}\right)$, where
$$
A(t,(\e^{i\theta},\xi))=\begin{cases}
  1+it^2u(\xi)u'(\xi)d\xi\,d\theta, & \text{if}\quad \xi\geq 0;\\
  1, & \text{if}\quad \xi\leq 0.
\end{cases}
$$

Thus $\eta(\sigma,\A,t)=- \str\left(iv_\sigma
\,\e^{\F(\sigma,\A,t)}\right)$ is given by
$$
 \eta(\sigma,\A,t)(\e^{i\theta},\xi)=\begin{cases}
  -2it\e^{-t^2 u(\xi)^2}u(\xi)^2 d\theta, & \text{if}\quad \xi\leq 0;\\
    0, & \text{if}\quad
   \xi\geq 0.
\end{cases}
$$
Finally, integrating  $\eta(\sigma,\A,t)$ in $t$ from $0$ to
$\infty$, we find that $\beta(\sigma,\A)$ (which is defined on
$\{(\e^{i\theta},\xi);\,u(\xi)\neq 0 \}= \T^*S^1\setminus   \supp(\sigma)$)  is equal to
\begin{equation}\label{eq:beta-TS}
\beta(\sigma,\A)(\e^{i\theta},\xi)=\begin{cases}
  -id\theta, & \text{if}\quad \xi\geq 0, u(\xi)\neq 0 \\
   0, & \text{if} \quad  \xi\leq 0, u(\xi)\neq 0 .
\end{cases}
\end{equation}

We have then proved the following

\begin{prop}
$\bullet$ The relative Chern class $\chr(\sigma)$ is represented
$(0,\beta(\sigma,\A))$.

$\bullet$ The Chern class with compact support $\chc(\sigma)$
is represented by the differential form $-i\,\mathbf{1}_{\geq 0}\, d\chi\wedge d\theta$
where $\chi\in\f(\Rbb)$ is compactly supported and equal to $1$ on $[-1,1]$, and
$\mathbf{1}_{\geq 0}$ is the characteristic function of the interval
$[0,\infty[$
\end{prop}

Note that the
differential form $-i\,\mathbf{1}_{\geq 0}\, d\chi\wedge d\theta$ is
of integral $-2i\pi$ on $\T^*S^1$ (which is oriented by $d\theta\wedge d\xi$).

%

\bigskip

\subsubsection{The space $\Rbb^2$.}\label{basic2}

Now we consider the case where $N=\Rbb^2\simeq\Cbb$. Take
$\Ecal^+=\Ecal^-$ the trivial bundles $N\times \Cbb$ over $N$. We
consider Bott's symbol $\sigma_b: \Ecal^+\to \Ecal^-$ which is given
by the map $\sigma_b(z)=z$ for $z\in N\simeq \Cbb$. The support of
$\sigma_b$ is reduced to the origin $\{0\}$, thus $\sigma_b$ defines
an element of $\KK^0(\Rbb^2)$. Recall that the Bott isomorphism
tells us that $\KK^0(\Rbb^2)$ is a free $\Zbb$-module with base
$\sigma_b$.

We choose on $\Ecal^\pm$ the trivial connections
$\nabla^+=\nabla^-=d$. Let   $\A=\nabla^+\oplus \nabla^-$  be the trivial connection  on
$\Ecal^+\oplus \Ecal^-$. The curvature $\F(\sigma_b,\A,t)$ of the
super-connection $\A^{\sigma_b}(t)= \left(
 \begin{array}{cc}
d & 0\\
0 & d\\
\end{array}\right)+\left(\begin{array}{cc}
0 & it\overline{z}\\
itz & 0\\
\end{array}\right)$
has the matrix form (see (\ref{eq:convention}))
$$
\F(\sigma_b,\A,t)=\left(\begin{array}{cc}-t^2|z|^2 & 0\\ 0 & -t^2|z|^2\\
\end{array}\right)-it \left(\begin{array}{cc} 0 & d\overline{z}\\
dz & 0\\ \end{array}\right).
$$

Thus
$$\e^{\F(\sigma_b,\A,t)}=\e^{-t^2|z|^2}\left(\begin{array}{cc}
1-\frac{t^2}{2}dz d\overline{z} & -itd\overline{z}\\
-itdz & 1+\frac{t^2}{2}dz d\overline{z}\\
\end{array}\right)$$
and $\eta(\sigma_b,\A, t)=- \str\left(iv_{\sigma_b}\,\e^{\F(\sigma_b,\A,t)}\right)$ is equal  to
%
%
\begin{equation}\label{eq:eta}
    \eta(\sigma_b,\A, t)=-t(\overline{z}dz-z d\overline{z})\e^{-t^2|z|^2}.
\end{equation}
When $z\neq 0$, we obtain that $\beta(\sigma_b,\A)(z)= \int_{0}^{\infty}\eta(\sigma_b,\A, t)dt$
is equal to $\frac{1}{2|z|^2}(z d\overline{z}-\overline{z}dz)=-i \,d (\arg z)$. Thus we have
\begin{equation}\label{chQb}
\chr(\sigma_b)=[0, -i\, d (\arg z)].
\end{equation}
It is easy to see that  $\chr(\sigma_b)$ is a basis of the vector space $\Hcal^*(\Cbb,\Cbb\setminus\{0\})$.

Take $f\in\f(\Rbb)$ with compact support and equal to $1$ in a
neighborhood of $0$.  Let $\chi(z):= f(|z|^2)$. Then the class
$\ch_c(\sigma_b)\in \Hcal^*_c(\Rbb^2)$ is represented by the
differential form $ c(\sigma_b,\A,\chi)=\chi\ch(\A)+
d\chi\beta(\sigma_b,\A)$. Here  the differential form $\ch(\A)$ is
identically equal to $0$. We obtain
\begin{eqnarray*}
  c(\sigma_b,\A,\chi) &=& d(f(|z|^2))\wedge\beta(\sigma_b,\A) \\
   &=& -f'(|z|^2) d\overline{z}\wedge dz.
\end{eqnarray*}

Remark that $c(\sigma_b,\A,\chi)$ is compactly supported and of
integral equal to $2i\pi$ on $\Rbb^2$ (with orientation $dx\wedge
dy$).  Thus $\frac{1}{2i\pi}c(\sigma_b,\A,\chi)$ is a
representative of the Thom form of $\Rbb^2$.

\begin{rema}
For $t>0$, the Chern character of the super-connection
$\A^{\sigma_b}(t)$ is  the degree $2$ differential form with
``Gaussian look''
$$\ch(\sigma_b,\A,t)=-\e^{-t^2|z|^2}t^2 dz d\overline{z}.$$

For any $t>0$,  $\ch(\sigma_b,\A,t)$ and $c(\sigma_b,\A,\chi)$
coincide in the cohomology  $\Hcal^*_{\dr}(\Rbb^2)$, as
follows from Proposition \ref{prop:chg=ch-dec-rapid}. In particular
they have the same integral.
\end{rema}

\subsubsection{The multiplicativity property on
$\Cbb^2$}\label{C2}

Following the notations of preceding example, we consider $\Cbb^2$
with coordinates $z=(z_1,z_2)$ and morphisms $\sigma_1=z_1$ and
$\sigma_2=z_2$. Then the tensor product morphism is
$$\sigma_1\odot \sigma_2=\left(%
\begin{array}{cc}
  z_1 & -\overline{z_2} \\
  z_2 & \overline{z_1} \\
\end{array}
\right).
$$

The morphism $\sigma_1\odot \sigma_2$ has support $z_1=z_2=0$.
A calculation similar to the calculation done in the preceding
section gives the following

\begin{prop}\label{basic3}
The relative chern class $\chr(\sigma_1\odot \sigma_2)\in
\Hcal^*(\Cbb^2,\Cbb^2\setminus(0,0))$
is represented by $(0, \beta_{12})$, where
$$
\beta_{12}=\frac{-1}{2|z|^4}
\Big((\overline{z_1}dz_1-z_1d\overline{z_1}) \wedge
d\overline{z_2}\wedge dz_2 +(\overline{z_2}dz_2-z_2d\overline{z_2})
\wedge d\overline{z_1}\wedge dz_1\Big)
$$
is a closed form on $\Cbb^2\setminus(0,0)$.
\end{prop}

Remark that $\beta_{12}$ is invariant under the symmetry group
$U(2)$ of $\Cbb^2$.

Recall that $\chr(\sigma_k)=[0,\beta_k]$, with
$\beta_k=-\frac{\overline{z_k}dz_k-z_kd\overline{z_k}}{2|z_k|^2}$.
The wedge product  $\beta_1\wedge \beta_2$ is not defined on
$\Cbb^2\setminus(0,0)$. Introduce a partition of unity
$\Phi_1,\Phi_2$ with respect to the covering $U_1\cup U_2$ of
$\Cbb^2\setminus(0,0)$, with $U_k=\{z,z_k\neq 0\}$. Then the
relative product $\chr(\sigma_1)\diamond \chr(\sigma_2)$ has
representative $(0,\beta)$, with $\beta=-d\Phi_1 \wedge \beta_1\wedge \beta_2$.

We now compute the forms $B_1,B_2$ of the equation
(\ref{eq:produit-relatif}). The form $\eta_k(t)$ have been
computed  (Equation \ref{eq:eta}). From this it is easy to compute
$B_1=\Phi_1\int_{0\leq s\leq t}\eta_1(t)\eta_2(s)ds\, dt$ and $B_2$.
We obtain

\begin{eqnarray*}
B_1&=&\Phi_1(z_1,z_2)
\frac{(\overline{z_1}dz_1-z_1d\overline{z_1})\wedge
(\overline{z_2}dz_2-z_2d\overline{z_2})}{4|z_1|^2(|z_1|^2+|z_2|^2)},\\
B_2&=&\Phi_2(z_1,z_2)
\frac{(\overline{z_1}dz_1-z_1d\overline{z_1})\wedge
(\overline{z_2}dz_2-z_2d\overline{z_2})}{4|z_2|^2(|z_1|^2+|z_2|^2)}.
\end{eqnarray*}

Here $B_1-B_2$ is a two form which is well defined on
$\Cbb^2\setminus(0,0)$ and the relation $\Phi_1+\Phi_2=1$ imply
$$
\beta_{12}-\beta=d(B_1-B_2).
$$
This shows that the class $\chr(\sigma_1\odot \sigma_2)$ is the
product  $[0,\beta_1]\diamond [0,\beta_2]$.

\medskip

We can now look at the different representatives of the Chern class with compact
support $\ch_c(\sigma_1\odot \sigma_2)\in\Hcal^*_c(\Cbb^2)$.
Let $f\in\f(\Rbb)$ with compact support and equal to $1$ in a
neighborhood of $0$. We consider the functions
$\chi(z)=f(|z_1|^2+|z_2|^2)$ and $\chi_k(z_k)=f(|z_k|^2)$. Let
$\Omega=d\overline{z_1} \wedge dz_1\wedge d\overline{z_2}\wedge
dz_2$.

\begin{prop}
The Chern class $\ch_c(\sigma_1\odot \sigma_2)\in\Hcal^*_c(\Cbb^2)$
is represented by any of the following differential forms
\begin{eqnarray*}
c(\sigma_1\odot\sigma_2,\A,\chi)&=& -\frac{f'\left(|z_1|^2+|z_2|^2\right)}{|z_1|^2+|z_2|^2}\, \Omega\\
c(\sigma_1,\A_1,\chi_1)\wedge c(\sigma_2,\A_2,\chi_2)&=& f'\left(|z_1|^2\right)
f'\left(|z_2|^2\right)\,\Omega.
\end{eqnarray*}

\end{prop}

Clearly the first representative is ``better'', as it is  invariant
by the full symmetry group $SO(4)$ on $\Cbb^2=\Rbb^4$.

\section{Riemann-Roch formula in relative cohomology}\label{sec:RR-formula}

In this section, we work with differential forms with {\em real} coefficients until Subsection \ref{notation2}.

\subsection{Some notations}\label{notation1}

Let $V$ be an Euclidean vector space of dimension $d$, with oriented orthonormal basis
$e_1,e_2,\ldots,e_d$.
We identify the Lie algebra $\sogot(V)$ of $SO(V)$ with $\Lambda^2 V$ as follows:
to an antisymmetric matrix $A$ in $\sogot(V)$, we associate  the element $\sum_{i<j}(Ae_i,e_j) e_i\wedge e_j$ 
of $\Lambda^2 V$. This identification will be in place throughout this section.

The Berezin integral $\bere: \Lambda V\to \Rbb$  is the $\Rbb$-linear map which vanishes on $\Lambda^i V$
for $i<d$ and is such that $\bere(e_1\wedge e_2\wedge \cdots\wedge e_d)=1$.

\subsection{Thom class in relative cohomology}\label{sec:Thom-class}

Let $M$ be a manifold. Let $p:\Vcal\to M$ be a real oriented
Euclidean vector bundle over $M$ of rank $d$. In this section, we give a
construction for the relative Thom form, analogous to Quillen's
construction of the Chern character. Here, we use the Berezin integral which
is the ``super-commutative" analog of the
super-trace for endomorphisms of a super-space.

Recall the sub-space $\Acal^*_{\cf}(\Vcal)\subset
\Acal^*(\Vcal)$ of (real) differential forms on $\Vcal$
which have a compact support in the fibers of $p:\Vcal\to M$. We
have also defined the sub-space $\Acal^*_{\dr}(\Vcal)$. The
integration over the fiber, that we denote by $p_*$,  is well defined on the three spaces
$\Acal^*(\Vcal,\Vcal\setminus M)$,
$\Acal^*_{\cf}(\Vcal)$ and
$\Acal^*_{\dr}(\Vcal)$ and take values in $\Acal^*(M)$.
A Thom form on $\Vcal$ will be
a (real) closed element which integrates to the constant function $1$ on
$M$.

Let $\nabla$ be an Euclidean  connection on $\Vcal$. As the structure group of $\Vcal$ is  
the Lie group $SO(V)$ with Lie algebra $\sogot(V)$,
the curvature $\F$  of $\nabla$ is   a two-form with values antisymmetric transformations of $\Vcal$. We
 will identify  the curvature to  an element
 $\F\in\Acal^2(M,\Lambda^2\Vcal)$    according to the isomorphism
 $\sogot(V)\sim \Lambda^2 V$ described above.
Let $\bere:\Gamma(M,\Lambda \Vcal)\to \f(M)$ be the Berezin integral that
we extend to a $\Rbb$-linear map $\bere:\Acal^*(M,\Lambda \Vcal)\to \Acal^*(M)$. The pfaffian of an element
$L\in \Acal^*(M,\Lambda^2 \Vcal)$ is defined by: $\Pf(L):=\bere(\e^L)$.

\begin{defi}
Let $\nabla$ be an  Euclidean  connection  on $\Vcal$, with curvature form
$\F$. The Euler form $\Eul(\Vcal,\nabla)\in \Acal^*(M)$ of the bundle  $\Vcal\to M$
is the closed real differential form on $M$ defined by
$\Eul(\Vcal,\nabla):=\Pf\left(-\frac{\F}{2\pi}\right)$.
The class of $\Eul(\Vcal,\nabla)$, which does not depend on $\nabla$, is denoted
by $\Eul(\Vcal)\in\Hcal^{d}(M)$.
\end{defi}

\begin{rema}
Since the pfaffian vanishes when the rank of $\Vcal$ is odd,
the Euler class $\Eul(\Vcal)\in\Hcal^{d}(M)$
is identically equal to $0$ when the rank of $\Vcal$ is odd.
\end{rema}

Let us consider the vector bundle $p^*\Vcal\to \Vcal$ equipped with the pull-back connection $p^*\nabla$. Let
$\x$ be the canonical section of the bundle $p^*\Vcal$. We consider the $\Zbb/2\Zbb$ graded algebra
$\Acal^*(\Vcal,\Lambda p^*\Vcal)$ which is equipped with the Berezin integral
$\bere:\Acal^*(\Vcal,\Lambda p^*\Vcal)\to\Acal^*(\Vcal)$.

Let $f_t^\nabla\in\Acal^*(\Vcal,\Lambda p^*\Vcal)$ be the element defined
by the equation
\begin{equation}\label{eq:f-t}
f_t^\nabla= -t^2\|\x\|^2 + t\, p^*\nabla\x+ \frac{1}{2} p^*\F.
\end{equation}

We consider the real differential forms on $\Vcal$
defined by
\begin{eqnarray}\label{eq:C-eta-wedge}
\mathrm{C}_{\wedge}^t&:=&\bere\left(\e^{f_t^\nabla}\right),\\
\eta_{\wedge}^t&:=& - \bere\left(\x\e^{f_t^\nabla}\right).
\end{eqnarray}
Here the exponentials are computed in the super-algebra
$\Acal^*(\Vcal,\Lambda p^*\Vcal)$.
To be more concrete, this calculation  is performed explicitly for a rank two bundle
in Example \ref{rank2} at the end of this subsection.

\begin{lemm}\label{thomgaussianclosed}
The differential form $ \mathrm{C}_{\wedge}^t$ is closed.
Furthermore,
\begin{equation}\label{eq:Thom}
\frac{d}{dt}\mathrm{C}_{\wedge}^t=-d(\eta_{\wedge}^t).
\end{equation}
\end{lemm}

\begin{proof}
The proof of the first point is given in \cite{BGV} (Chapter 7,
Theorem 7.41). We recall the proof. We denote by
$\iota_{\wedge}(\x)$ the derivation of the super-algebra $\Acal^*(\Vcal,\Lambda p^*\Vcal)$ such that
$\iota_{\wedge}(\x)s=\langle \x,s\rangle$ when $s\in \Acal^0(\Vcal,\Lambda^1 p^*\Vcal)$.
We extend the connection $p^*\nabla$ to a derivation $\nabla^\wedge$ of $\Acal^*(\Vcal,\Lambda p^*\Vcal)$.
We consider the derivation $\nabla^\wedge -2t \iota_{\wedge}(\x)$ on $\Acal^*(\Vcal,\Lambda p^*\Vcal)$.
It is easy to verify that
\begin{equation}\label{eft}
(\nabla^\wedge -2t \iota_{\wedge}(\x))f_t^\nabla=0.
\end{equation}
Then, the exponential $\e^{f_t^\nabla}$ satisfies also
$(\nabla^\wedge-2t\iota_{\wedge}(\x))(\e^{f_t^\nabla})=0$. The Berezin integral is such
that $\bere(\iota_{\wedge}(\x)\alpha)=0$ and $\bere(\nabla^\wedge\alpha)=d(\bere(\alpha))$ for any $\alpha\in
\Acal^*(\Vcal,\Lambda p^*\Vcal)$. This shows that $d\left(\bere(\e^{f_t^\nabla})\right)=0$.

Let us prove the second point. We have
$d\circ \bere \big(\x\e^{f_t^\nabla}\big)=$ \break
$\bere\circ\, (\nabla^\wedge-2t\iota_{\wedge}(\x))\big(\x\e^{f_t^\nabla}\big)$, and
since $(\nabla^\wedge -2t\iota_{\wedge}(\x))\e^{f_t^\nabla}=0$, we get
\begin{eqnarray*}
 (\nabla^\wedge-2t\iota_{\wedge}(\x))\big(\x\e^{f_t^\nabla}\big)&=&
 \big((\nabla^\wedge-2t \iota_{\wedge}(\x))\cdot \x\big)\e^{f_t^\nabla}\\
&=& (\nabla^\wedge \x-2t\|\x\|^2)\e^{f_t^\nabla}\\
&=&\frac{d}{dt}\e^{f_t^\nabla}.
\end{eqnarray*}
\end{proof}

When $t=0$, then $\mathrm{C}_{\wedge}^0$ is just equal to
$\Pf(\frac{\F}{2})=(-\pi)^{d/2}\Eul(\Vcal,\nabla).$ When $t=1$, then
$\mathrm{C}_{\wedge}^1=\bere(\e^{f_1^\wedge})=\e^{-\|\x\|^2}Q$ is a
closed form with a Gaussian look on $\Vcal$ :
$Q$ a differential form on $\Vcal$  with a polynomial growth on the fiber of $\Vcal\to M$ (we will be more
explicit in a short while). This differential form was considered by
Mathai-Quillen in \cite{Mathai-Quillen}.


\bigskip

We have $\eta_{\wedge}^t=\e^{-t^2\|\x\|^2}Q(t)$ where $Q(t)$ is a differential form on $\Vcal$
with a polynomial growth on the fiber of $\Vcal$ and which depends polynomially on $t\in \Rbb$.
Thus, if $\x\neq 0$, when $t$ goes to infinity, $\eta_{\wedge}^t$ is an exponentially decreasing
function of $t$. We can thus define the following differential form
on $\Vcal\setminus M$ :
\begin{equation}\label{def:beta-form}
\beta_{\wedge}=\int_{0}^{\infty}\eta_{\wedge}^t\,dt.
\end{equation}
If we integrate (\ref{eq:Thom}) between  $0$ and $\infty$, we get
$\mathrm{C}_{\wedge}^0=d(\beta_{\wedge})$ on $\Vcal\setminus M$.
Thus the couple $(\mathrm{C}^0_\wedge,\beta_\wedge)$ defines a canonical
relative class
\begin{equation}\label{def:pf.beta}
\left[\Pf(\hbox{$\frac{\F}{2}$}),\beta_\wedge\right]\in\Hcal^*(\Vcal, \Vcal\setminus M)
\end{equation}
of degree equal to the rank of $\Vcal$.

We give the explicit formula for this relative class in the case  of  a rank two 
Euclidean bundle in Example \ref{rank2}.

Consider now the cohomology with compact support in the fiber of
$\Vcal$. Let $\mathrm{C}_\Vcal$ be the image of $\left[\Pf(\hbox{$\frac{\F}{2}$}),\beta_\wedge\right]$
through the map $\p_{\cf}:\Hcal^*(\Vcal, \Vcal\setminus M)\to \Hcal^*_{\cf}(\Vcal)$.

\begin{prop}\label{prop:Thom}
Let $\chi\in\f(\Vcal)$ be a function with
compact support in the fibers and equal to $1$ in a neighborhood of $M$. The
form
$$
\mathrm{C}_\Vcal^\chi= \chi \Pf(\hbox{$\frac{\F}{2}$})
 +d\chi\beta_{\wedge}
$$
is a closed differential form with compact support in the fibers on $\Vcal$.  Its cohomology class
in $\Hcal^*_{\cf}(\Vcal)$ coincides with
$\mathrm{C}_\Vcal$ : in particular, it does
not depend on the choice of $\chi$. We have $\frac{1}{\epsilon_d}p_*\left(\mathrm{C}^\chi_\Vcal\right)=1$,
with $\epsilon_d=(-1)^{\frac{d(d-1)}{2}}\pi^{d/2}$. Thus
$\frac{1}{\epsilon_d}\mathrm{C}^\chi_\Vcal$ is a Thom form in $\Acal^*_{\cf}(\Vcal)$.
\end{prop}
\begin{proof}
The first assertions are consequence of the definition of
$\mathrm{C}_\Vcal$. To compute $p_*\left(\mathrm{C}^\chi_\Vcal\right)$, we
may choose  $\chi= f(\|\x\|^2)$ where $f\in \f(\Rbb)$ has a compact
support and is equal to $1$ in a neighborhood of $0$. We work with a local oriented orthonormal frame
 $(e_1,\ldots,e_d)$ of $\Vcal$ : we have
$\x=\sum_i x_i e_i$, and $p^*\nabla\x=\sum_i dx_i e_i+x_i p^*\nabla e_i$.

The component of maximal degree in the fibers of the differential form $\eta_\wedge^t$ is
$(-1)^{\frac{d(d-1)}{2}}t^{d-1} \e^{-t^2\|x\|^2} \sum_k (-1)^k x_k dx_{1}
\cdots \widehat{dx_k}\cdots dx_{d}$ (see Proposition \ref{prop:explicit}). Then, the
component of maximal degree in the fibers of the differential form $d\chi\wedge
\eta^t_\wedge$ is
$$
-2(-1)^{\frac{d(d-1)}{2}}t^{d-1}f'(\|x\|^2)\|x\|^2 \e^{-t^2\|x\|^2} dx_{1}\cdots dx_{d}.
$$
Hence, using the change of variables $x\to
\frac{1}{t}x$,
\begin{eqnarray*}
  p_*\left(\mathrm{C}_\Vcal^\chi\right) &=&
  -2(-1)^{\frac{d(d-1)}{2}}\int_{0}^\infty t^{d-1}\left(\int_{\Rbb^d}
  f'(\|x\|^2)\|x\|^2 \e^{-t^2\|x\|^2} dx\right)dt\\
   &=&(-1)^{\frac{d(d-1)}{2}}\int_{0}^\infty \underbrace{\left(\int_{\Rbb^d}
     f'(\frac{\|x\|^2}{t^2})(\frac{-2\|x\|^2}{t^3})\e^{-\|x\|^2}
  dx\right)}_{I(t)}dt.
\end{eqnarray*}
Since for $t>0$, $I(t)=\frac{d}{dt}\left(\int_{\Rbb^d}
  f(\frac{\|x\|^2}{t^2})\e^{-\|x\|^2} dx\right)$, we have
$p_*\left(\mathrm{C}_\Vcal^\chi\right)=$ \break $(-1)^{\frac{d(d-1)}{2}}\int_{\Rbb^d} \e^{-\|x\|^2}
dx= (-1)^{\frac{d(d-1)}{2}}\pi^{d/2}.$
\end{proof}

Using the differential form $\mathrm{C}_{\wedge}^1$, it is possible to construct
representatives of a Thom form with Gaussian look.

\begin{prop}[Mathai-Quillen]\label{thomgaussian}
The differential form $\mathrm{C}_{\wedge}^1$ is a closed  form which belongs to
$\Acal_{\dr}^*(\Vcal)$. We have $\frac{1}{\epsilon_d}p_*\left(\mathrm{C}_{\wedge}^1\right)=1$,
with $\epsilon_d=(-1)^{\frac{d(d-1)}{2}}\pi^{d/2}$. Thus
$\frac{1}{\epsilon_d}\mathrm{C}_{\wedge}^1$ is a Thom form in
$\Acal_{\dr}^*(\Vcal)$.
\end{prop}

\begin{proof}
 By Proposition \ref{prop:explicit},  the component of maximal degree in the fibers of the differential
form $\mathrm{C}_{\wedge}^1$ is the term 
$(-1)^{\frac{d(d-1)}{2}} \e^{-\|x\|^2}dx_1 \cdots dx_d$.
\end{proof}

We summarize Propositions \ref{prop:Thom} and \ref{thomgaussian}
in the following theorem.

\begin{theo}\label{theo:thom-form}
Let $p:\Vcal\to M$ be an oriented Euclidean vector bundle of rank $d$ equipped with an Euclidean
connection $\nabla$, with curvature $\F$. Let
$\epsilon_d:=(-1)^{\frac{d(d-1)}{2}}\pi^{d/2}$.
Let $\bere:\Acal^*(\Vcal,\Lambda p^*\Vcal)\to \Acal^*(\Vcal)$
be the Berezin integral. Let
\begin{eqnarray*}
f_t^\nabla&=& -t^2\|\x\|^2 + tp^*\nabla \x + \frac{1}{2}p^*\F,\\
\eta_{\wedge}^t&=&- \bere\left(\x\e^{f_t^\nabla}\right),\\
\beta_{\wedge}&=&\int_{0}^{\infty}\eta_{\wedge}^tdt.
\end{eqnarray*}

$\bullet$ $\tur(\Vcal,\nabla)=\frac{1}{\epsilon_d}
\left(\Pf(\hbox{$\frac{\F}{2}$}),\beta_{\wedge}\right)$
is a Thom form in $\Acal^*(\Vcal,\Vcal\setminus M)$. It defines a
Thom class
$$
\tur(\Vcal)\in\Hcal^*(\Vcal,\Vcal\setminus M).
$$

$\bullet$ $\tuc(\Vcal,\nabla,\chi)=\frac{1}{\epsilon_d}\mathrm{C}_{\Vcal}^\chi=
\frac{1}{\epsilon_d}\left( \chi \Pf(\hbox{$\frac{\F}{2}$}) +d\chi \beta_\wedge\right)$
 is a  Thom form in $\Acal_{c}^*(\Vcal)$. Here $\chi\in\f(\Vcal)$ is a function with compact support
 in the fibers of $\Vcal$  and equal to $1$ in a neighborhood of $M$. It defines a Thom class
$$
\tuc(\Vcal)\in\Hcal^*_c(\Vcal).
$$

$\bullet$
The Mathai-Quillen form $\tumq(\Vcal,\nabla)=\frac{1}{\epsilon_d}\mathrm{C}_{\wedge}^1=
\frac{1}{\epsilon_d}\bere\left(\e^{f_1^\nabla}\right)$
is a Thom  form in $\Acal_{\dr}^*(\Vcal)$.
It defines a Thom class
$$\tumq(\Vcal)\in\Hcal^*_{\dr}(\Vcal).$$

\end{theo}

Thus the use of the Berezin integral allowed us to give slim
formulae for Thom forms in relative cohomology, as well as in
compactly supported cohomology or in rapidly decreasing
cohomology.

\begin{exa}\label{rank2}{\rm  Vector bundle of rank $2$}.

We write explicitly the formulae of this subsection  in  the case of an  Euclidean bundle $\Vcal\to M$  of rank $2$ in a local frame.
 Let  $(e_1,e_2)$ be a local oriented orthonormal frame. Let $\nabla$ be an Euclidean connection
 on $\Vcal$, so that $\nabla e_1=\eta e_2$, $\nabla e_2=-\eta e_1$,
               where $\eta$ is a real valued one form on $M$.
Then $$\F= d\eta (e_1\wedge e_2),\hspace{1cm} p^*\nabla \x=\eta_1 e_1+\eta_2 e_2$$
with $\eta_1=dx_1-x_2 \eta, \eta_2=dx_2 +x_1 \eta$, and
$$f_t^\nabla= -t^2\|x\|^2 + t(\eta_1 e_1+\eta_2 e_2)+\frac{1}{2} d\eta (e_1\wedge e_2).$$

The exponential of $f_t^{\nabla}$ in the super-algebra
$\Acal^*(\Vcal,\Lambda p^*\Vcal)$
is

$$\e^{f_t^{\nabla}}=\e^{-t^2\| x\|^2}\left(1+ \frac{d\eta}{2} e_1\wedge e_2+ t(\eta_1 e_1+\eta_2 e_2) -t^2(\eta_1\wedge \eta_2) e_1\wedge e_2\right).$$

Thus we have the formulae:

\begin{eqnarray*}
\mathrm{C}_{\wedge}^t&=&\e^{-t^2\| x\|^2}
(\frac{d\eta}{2}-t^2 \eta_1\wedge \eta_2),\\
\eta_{\wedge}^t&=& t\e^{-t^2\| x\|^2} (x_1 \eta_2-x_2\eta_1),\\
\beta_{\wedge}&=&\frac{x_1\eta_2-x_2 \eta_1}{2\| x\|^2}.
\end{eqnarray*}

So Thom forms  are given by

$$
\bullet\quad
\tur(\Vcal)=\frac{-1}{2\pi}\left [d\eta,\eta + \frac{x_1 dx_2-x_2 dx_1}{\| x\|^2}\right],
$$

$$
\bullet\quad
\tuc(\Vcal)=\frac{-1}{2\pi}\Big( 2f'(\| x\|^2) dx_1\wedge dx_2  + d\left(f(\| x\|^2) \wedge \eta\right)\Big),
$$
where $f$ is a compactly supported function on $\Rbb$, identically equal to $1$ in a neighborhood of
$0$,

$$
\bullet\quad
\tumq(\Vcal)=\frac{1}{2\pi}\e^{-\| x\|^2}\Big( 2 dx_1\wedge dx_2 - d\eta +
d\left(\| x\|^2\right) \wedge \eta \Big).
$$

\end{exa}

\medskip

\subsection{Explicit formulae for the Thom forms of a vector bundle}\label{subsec:explicit-Thom-vector}

Let us give  explicit local formulae  for a general Euclidean vector bundle.

Given a local oriented orthonormal frame $(e_1,\ldots,e_d)$ of the vector bundle $\Vcal$, we work with the identification
$(m,x)\mapsto \sum_i x_i e_i(m)$ from $M\times \Rbb^d$ into $\Vcal$. The element $p^*\nabla\x$ is then
equal to $\sum_i \eta_i e_i$ with $\eta_i=dx_i + \sum_k x_k(\nabla e_k, e_i)$.

If $I=[i_1,i_2,\ldots,i_p]$ (with $i_1<i_2<\cdots<i_p$) is a
subset of $[1,2,\ldots,d]$, we use the notations $e_I=e_{i_1}\wedge \cdots \wedge e_{i_p}$ and
$\eta_I=\eta_{i_1}\wedge \cdots \wedge \eta_{i_p}$. The curvature $\F$ decomposes as $\F:=\sum_{i<j} \F_{ij} e_i\wedge e_j$.
For any subset $I$ of  $[1,2,\ldots,d]$, we consider the two form $\F_I:=\sum_{i<j,i\in I,j\in I} \F_{ij} e_i\wedge e_j$ with values in
$\Lambda \Vcal_I$, where $\Vcal_I$ is the sub-bundle generated by the $e_i,i\in I$. Let $\Pf(\F_I)$ be its pfaffian. One sees easily that
\begin{equation}\label{eq:pfaffian-gene}
\e^{\frac{\F}{2}}=\sum_{I}\Pf(\hbox{$\frac{\F_I}{2}$})\, e_I \quad {\rm in}\quad \Acal^*(M,\Lambda \Vcal).
\end{equation}
Only those $I$ with $|I|$ even will
contribute to the sum (\ref{eq:pfaffian-gene}),  as otherwise the pfaffian of $\F_I$
vanishes.

If $I$ and $J$ are two disjoint subsets of $\{1,2,\ldots,d\}$, we denote by
$\epsilon(I,J)$ the sign such that $e_I\wedge e_J=\epsilon(I,J)\, e_{I\cup J}$.

\begin{prop}\label{prop:explicit}

$\bullet$ We have
$\tur(\Vcal)=\frac{1}{\epsilon_d}\left[\Pf(\hbox{$\frac{\F}{2}$}),\beta_{\wedge}\right]$
with
$$
\beta_\wedge=\sum_{k,I,J} \gamma_{(k,I,J)}
 \Pf\left(\hbox{$\frac{\F_I}{2}$}\right)\frac{x_k \eta_J}{\|x\|^{|J|+1}},
$$
with
$$
\gamma_{(k,I,J)}=-\frac{1}{2}(-1)^{\frac{|J|(|J|+1)}{2}}\Gamma\left(\hbox{$\frac{|J|+1}{2}$}\right)
\epsilon(I,J)\epsilon(\{k\},I\cup J).
$$
Here for $1\leq k\leq d$, the sets $I,J$ vary over the subsets of
$\{1,2,\ldots,d\}$ such that  $\{k\}\cup I\cup
J$ is a partition of $\{1,2,\ldots,d\}$. Only those $I$ with $|I|$ even will
contribute to the sum.

$\bullet$ The class $\tuc(\Vcal)$ is represented by the closed differential form
$$
\frac{1}{\epsilon_d}\Big(f(\|x\|^2) \Pf(\hbox{$\frac{\F}{2}$})+ 2f'(\|x\|^2) (\sum
x_i dx_i)\beta_{\wedge}\Big)
$$
where $f$ is a compactly supported function on $\Rbb$, identically equal to $1$ in a neighborhood of
$0$.

$\bullet$ We have
$$
\tumq(\Vcal)=\frac{1}{(\pi)^{d/2}}\e^{-\|x\|^2} \sum_{I} (-1)^{\frac{| I |}{2}}\epsilon(I,I')
\Pf(\hbox{$\frac{\F_I}{2}$}) \eta_{I'}.
$$
Here $I$ runs over the subset of $\{1,2,\ldots,d\}$ with an even number of elements,
and $I'$ denotes the complement of $I$.
\end{prop}

\begin{proof}

It follows from the explicit description  of our forms  and from
the  formula $\int_0^{\infty} \e^{-t^2} t^adt=\frac{1}{2}
\Gamma(\frac{a+1}{2})$.

\end{proof}

\subsection{More notations}\label{notation2}
We recall notations from \cite{BGV}.
Let $V$ be an Euclidean vector space of even dimension $d=2n$ with oriented orthonormal basis
$e_1,e_2,\ldots,e_d$.
Let $C(V)$ be the Clifford algebra of $V$. Then $C(V)$ is generated by elements $c_i$ with relations $c_ic_j+c_jc_i=0$, for $i\neq j$, and $c_i^2=-1$.
We denote by $\Sigma:C(V)\to \Lambda V$ the symbol isomorphism. Thus, for
$1\leq i_1<i_2<\cdots<i_k\leq d$, we have  $\Sigma(c_{i_1}c_{i_2}\cdots c_{i_k})=
e_{i_1}\wedge e_{i_2}\wedge \cdots \wedge e_{i_k}$. Let $C^{[i]}(V)=\Sigma^{-1}(\Lambda^iV)$.
We denote by $\tau: C^{[2]}(V)\to \sogot(V)$ the map such that $\tau(c)v=cv-vc$, for
$c\in  C^{[2]}(V)$ and $v\in V$. Then $\tau(c_i c_j)(e_i)=2 e_j$, for $i\neq j$.
We denote by $S=S^+\oplus S^-$ the complex spinor space. We denote by $\cbf$ the Clifford action of $C(V)$ on $S$.  If $v\in V$, then $\cbf(v)$ on $S$ interchanges $S^+$ and $S^-$ and satisfies $\cbf(v)^2=-\|v\|^2{\rm Id}_S$. The supertrace of the action of the even element $c_1c_2\ldots c_d$ on $S$ is $(-2i)^n$.

Let $V=\Rbb e_1\oplus \Rbb e_2$ be of dimension $2$ .
  We consider the super-algebra $A\otimes C(V)$
where $A$ is a super-commutative algebra.
Then for $a_1, a_2$ odd  elements in $A$,
 and $b$ an even element of $A$,
we have

$$\exp( a_1 c_1+a_2 c_2+ b c_1 c_2)=\cos(b)+\sin(b) c_1 c_2\hspace{3cm}$$
$$+\frac{\sin(b)}{b}  (a_1 c_1+a_2c_2)+
\frac{\sin(b)-b \cos(b)}{b^2}  a_1 a_2-\frac{\sin(b)}{b} a_1 a_2 c_1 c_2.$$

This formula can be verified using, for example, the differential equation $\frac{d}{dt} \exp(tX)=X\exp(tX)$  for the exponential (see also \cite{BGV}, proof of  Proposition 7.43).

\subsection{Riemann-Roch relation}\label{sec:Riemman-Roch}

Let $p:\Vcal\to M$ be an oriented Euclidean vector bundle of even rank $d=2n$ with spin structure, and let $\Scal\to M$ be the corresponding
spin super-bundle. Let $C(\Vcal)\to M$ be the  Clifford bundle. We denote by
$$
\cbf : C(\Vcal)\to \End_\Cbb(\Scal)
$$
the bundle map defined by the spinor representation.

 The vector
bundle $\Scal$ is provided with an Hermitian metric such that
$\cbf(v)^*=-\cbf(v)$ for $v\in\Vcal$. Consider the morphism
$\sigma_\Vcal:p^*\Scal^+\to p^*\Scal^-$ defined by
$$
\sigma_\Vcal:= -i\cbf(\x)
$$
where $\x: \Vcal\to p^*\Vcal$ is the canonical section. Then the odd linear map $v_{\sigma_\Vcal}:p^*\Scal\to
p^*\Scal$ is equal to $-i\cbf(\x)$.

Let $\nabla$ be an Euclidean  connection on $\Vcal$.  Then as explained in
Subsection \ref{sec:Thom-class}, the curvature $\F$ of $\nabla$ may be
 identified to an element of $\Acal(M,\Lambda^2\Vcal)$.
Thus $\Sigma^{-1}\F$ is an element of  $\Acal(M,C(\Vcal))$,
where $\Sigma: C(\Vcal)\to \Lambda \Vcal$ is  the symbol bundle map.
The connection $\nabla$  induces a connection $\nabla^\Scal$ on $\Scal$ with curvature $\F^{\Scal}$ (see Lemma \ref{lem:F-Scal}).
We work with the family of super-connections on $p^*\Scal$:
$$
\A^{\sigma}_t:=p^*\nabla^\Scal + t\cbf(\x)
$$
We see that the curvature of the super-connection $\A^{\sigma}_t$ is the even element $\F_t^{\cbf}\in \Acal^*(\Vcal, p^*\End_\Cbb(\Scal))$,
given by
$$
\F_t^{\cbf}=- t^2\|\x\|^2 +t\cbf(p^*\nabla\x) +p^*\F^\Scal
$$
where $\F^\Scal\in \Acal^2(M, \End_\Cbb(\Scal))$ is the curvature of $\nabla^\Scal$.

\begin{lemm}\label{lem:F-Scal}
$\bullet$ The following relation holds in $\Acal^2(M, \End_\Cbb(\Scal))$:
$$
\F^\Scal=\frac{1}{2}\cbf(\Sigma^{-1}\F).
$$
$\bullet$ We have $\F_t^{\cbf}=\cbf(\tilde{f}_t)$, where $\tilde{f}_t\in \Acal^*(\Vcal, p^*C(\Vcal))$ is given by
$$
\tilde{f}_t=- t^2\|\x\|^2 +t p^*\nabla\x +\frac{1}{2}p^*\Sigma^{-1}\F.
$$
$\bullet$ The image of $\tilde{f}_t$ by the bundle map $\Sigma$ is equal to the map
$f^\nabla_t\in \Acal^*(\Vcal, p^*\Lambda\Vcal)$ defined in (\ref{eq:f-t}).
\end{lemm}

We consider now in parallel the closed differential forms
$$
\ch(\sigma_\Vcal,\nabla^\Vcal,t):=\str\left(\expo(\F_t^{\cbf})\right),\quad
\mathrm{C}_\wedge^t:=\bere\left(\expo_{\wedge}(f^\nabla_t)\right).
$$
In the first case the exponential is computed in the super-algebra \break
$\Acal^*(\Vcal,p^*\End_\Cbb(\Scal))$, and the forms $\ch(\sigma_\Vcal,\nabla^\Vcal,t)$ have {\em complex} coefficients. In the second case, the exponential is computed in the super-algebra
$\Acal^*(\Vcal,p^*\Lambda\Vcal)$, and the forms $\mathrm{C}_\wedge^t$ have {\em real} coefficients.

In Example \ref{rank2cliff}, we will perform the explicit calculation of $\expo(\F_t^{\cbf})$    for a bundle of rank two.

We also consider in parallel the differential forms
$$
\eta^t_\cbf:=-\str\left(\cbf(\x)\expo(\F_t^{\cbf})\right),\quad
\eta_\wedge^t:=-\bere\left(\x\cdot\expo_{\wedge}(f^\nabla_t)\right).
$$
Note that the forms $\eta^t_\cbf$ have {\em complex} coefficients, and that
the forms $\eta_\wedge^t$ have {\em real} coefficients.

In the next definition, we return to the original definition of the curvature $\F$  of the Euclidean connection  $\nabla$, that is we consider $\F$  as a $2$-form with values antisymmetric transformations of $\Vcal$.
\medskip

\begin{defi}
We associate to the vector bundle $\Vcal$, equipped with the connection $\nabla$, the closed
real differential form on $M$ defined by
$$
\Ahat(\nabla):=\det{}^{1/2}
\left(\dfrac{\F}{\e^{\frac{\F}{2}}-\e^{-\frac{\F}{2}}}\right).
$$
Its cohomology class $\Ahat(\Vcal)\in\Hcal^{*}(M)$ is the $\Ahat$-genus of $\Vcal$.
\end{defi}

\begin{prop}
We have the following equalities of differential forms on $\Vcal$:
$$
\ch(\sigma_\Vcal,\nabla^\Vcal,t)=(-2i)^{n} ~ \Ahat(\nabla)^{-1} ~\mathrm{C}_\wedge^t
$$
and
$$
\eta^t_\cbf=(-2i)^{n}~ \Ahat(\nabla)^{-1} ~
\eta_\wedge^t.
$$
\end{prop}

\begin{proof} The proof of the first relation is done in \cite{BGV}, Section 7.
The same proof works for the second equality.
Let us give here a brief idea of the proof. Let $\str_C$ be the super-trace on $\Acal^*(\Vcal,p^*C(\Vcal))$
such that $\str_C(a)=\str(\cbf(a))$ for any element $a\in \Acal^*(\Vcal,p^*C(\Vcal))$. We have then to show that
\begin{equation}\label{eq:lem-rel-1}
\str_C\left(\expo_C(\tilde{f}_t)\right)=(-2i)^{n} ~ \Ahat(\nabla)^{-1} ~\bere\left(\expo_{\wedge}(\Sigma\tilde{f}_t)\right)
\end{equation}
and
\begin{equation}\label{eq:lem-rel-2}
\str_C\left(\x\cdot\expo_C(\tilde{f}_t)\right)= (-2i)^{n} ~ \Ahat(\nabla)^{-1} ~
\bere\left(\x\cdot\expo_{\wedge}(\Sigma\tilde{f}_t)\right).
\end{equation}

If $V$ is an oriented Euclidean vector space of even dimension $2n$, we have the following fundamental relation between
$\str_C(\expo_C(a))$ and $\bere(\expo_\wedge(\Sigma a))$  for $a\in C^2(V)$:

\begin{equation}\label{eq:lem-rel-3}
\str_C(\expo_C(a))=(-2i)^{n} \det{}^{1/2}
\left(\dfrac{\e^{\frac{\tau(a)}{2}}-\e^{-\frac{\tau(a)}{2}}}{\tau(a)}\right)\bere(\expo_\wedge(\Sigma a))
\end{equation}
(see \cite{BGV}, Section 3). We see then that (\ref{eq:lem-rel-1}) is an extension of (\ref{eq:lem-rel-3}) to the case
where $a\in A^-\otimes C^1(V)+A^+\otimes C^2(V)$ (here $A$ is a super-commutative super-algebra).  This is verified by an explicit computation when $V$ is of dimension $2$, using the formula for the exponential that we  recalled in  Subsection \ref{notation2}.

\end{proof}

We can now conclude with the

\begin{theo}\label{prop:chern-V-spin}
$\bullet$ We have the following equality in
$\Hcal^*(\Vcal, \Vcal\setminus M)$:
\begin{equation}\label{eq:chernr-V}
\chr(\sigma_\Vcal) =(2i\pi)^{n} ~ \Ahat(\Vcal)^{-1} ~ \tur(\Vcal).
\end{equation}

$\bullet$ We have the following equality in
$\Hcal^*_c(\Vcal)$:
\begin{equation}\label{eq:chernc-V}
\chc(\sigma_\Vcal) =(2i\pi)^n ~ \Ahat(\Vcal)^{-1} ~\tuc(\Vcal).
\end{equation}

$\bullet$ We have the following equality in
$\Hcal^*_{\dr}(\Vcal)$:
\begin{equation}\label{eq:cherndr-V}
\ch_Q(\sigma_\Vcal) =(2i\pi)^n ~ \Ahat(\Vcal)^{-1} ~\tumq(\Vcal).
\end{equation}
\end{theo}

Remark that these identities holds at the level of the
representatives.

\begin{exa}\label{rank2cliff} {\rm Vector bundle of rank two.}

We return to Example \ref{rank2}, and keep the same notations.
Then we have
$$\tilde f_t=-t^2\|x\|^2+ t (\eta_1 c_1+\eta_2 c_2)+\frac{1}{2} (d \eta)  c_1 c_2.$$

 We  use the formula for the exponential recalled in Subsection \ref{notation2}, and we obtain

\begin{eqnarray*}
\lefteqn{
\e^{\tilde f_t}=\e^{-t^2\|x\|^2}\Big(\cos (\frac{d \eta}{2})+\sin(\frac{d \eta}{2}) c_1 c_2+t
\frac{\sin(\frac{d \eta}{2})}{(\frac{d \eta}{2})}(\eta_1 c_1+\eta_2 c_2)
} \nonumber\\
& &+\  t^2
\frac{\sin(\frac{d \eta}{2})-(\frac{d \eta}{2}) \cos(\frac{d \eta}{2})}{(\frac{d \eta}{2})^2}  \eta_1 \eta_2 -t^2\frac{\sin(\frac{d \eta}{2})}{(\frac{d \eta}{2})} \eta_1 \eta_2 c_1 c_2\Big).
\end{eqnarray*}

The supertrace of the action of $c_1c_2$ on $S$ is $-2i$.
Thus we obtain
\begin{eqnarray*}
\ch(\sigma_\Vcal,\nabla^\Vcal,t)&=&(-2i)
\e^{-t^2\|x\|^2}\Big(\sin(\frac{d \eta}{2}) -t^2\frac{\sin(\frac{d \eta}{2})}{(\frac{d \eta}{2})} \eta_1 \eta_2 \Big),\\
\eta^t_\cbf&=&(-2i) \ t \e^{-t^2\|x\|^2} \frac{\sin(\frac{d \eta}{2})}{(\frac{d \eta}{2})} \Big(x_1 \eta_2-x_2\eta_1 \Big ),\\
\beta_{\cbf}&=&(-2i) \frac{\sin(\frac{d \eta}{2})}{(\frac{d \eta}{2})} \frac{(x_1 \eta_2-x_2\eta_1)}{2\|x\|^2}
\end{eqnarray*}

Finally, the relative  Chern character form  associated to $\sigma_\Vcal$ is

\begin{equation}
\chr(\sigma_\Vcal) =(-i) \frac{\sin(\frac{d \eta}{2})}{(\frac{d \eta}{2})}
 \left[d\eta, \eta +\frac{(x_1 dx_2-x_2dx_1)}{\|x\|^2}\right].
\end{equation}

We have $$\det{}^{1/2}
\left(\dfrac{\F}{\e^{\frac{\F}{2}}-\e^{-\frac{\F}{2}}}\right)=  \frac{(\frac{d \eta}{2})}
{\sin(\frac{d \eta}{2})}.$$
Thus we see that we have the relations
\begin{equation}
\chr(\sigma_\Vcal) =(2i\pi) ~ \Ahat(\Vcal)^{-1} ~\tur(\Vcal),
\end{equation}

\begin{equation}
\chc(\sigma_\Vcal) =(2i\pi) ~ \Ahat(\Vcal)^{-1} ~\tuc(\Vcal),
\end{equation}

and
\begin{equation}
\ch_Q(\sigma_\Vcal) =(2i\pi) ~ \Ahat(\Vcal)^{-1} ~\tumq(\Vcal)
\end{equation}

\end{exa}

at the level of differential forms.

\section{Appendix}

We give  a proof of the estimate used in this article. It is based
on Volterra's expansion formula: if $H$ and $R$ are elements in a
finite dimensional associative algebra, then
$\e^{(H+R)}=\e^{H}+\sum_{k=1}^{\infty}I_k(H,R)$ where
\begin{equation}\label{volterra}
I_k(H,R)= \int_{\Delta_k}\e^{s_1H}R\e^{s_2H}R\cdots
R\e^{s_{k}H}R\e^{s_{k+1}H}ds_1\cdots ds_{k}
\end{equation}
Here $\Delta_{k}$ is the simplex $\{s_i\geq 0\,;\,
s_1+s_2+\cdots+s_{k}+s_{k+1}=1\}$ which has the volume $\frac{1}{k!}$
for the measure $ds_1\cdots ds_{k}$.

Now, let $A=\oplus_{i=0}^R A_i$  be a complex finite dimensional
graded commutative algebra with a norm $\|\cdot\|$ such that
$\|ab\|\leq \|a\|\|b\|$. We denote by $A_+=\oplus_{i=1}^RA_i$.
Thus $\omega^{R+1}=0$ for any $\omega \in A_+.$ Let $E$ be a
finite dimensional  Hermitian vector space. Then $\End(E)\otimes
A$ is an algebra with a norm still denoted by $\|\cdot\|$. If
$H\in \End(E)$, we denote also by $H$ the element $H\otimes 1$ in
$\End(E)\otimes A$.

We denote by $\herm(E)\subset\End(E)$ the subspace formed by the
Hermitian endomorphisms. When $H\in\herm(E)$, we denote by
$\sm(H)\in\Rbb$ the smallest eigenvalue of $H$ : we have
$$
\Big|\!\Big| \e^{-H}\Big|\!\Big|=\e^{-\sm(H)},\quad \mathrm{for\
all} \quad H\in\herm(E).
$$

\begin{lemm}\label{suffit}
Let $\Pcal(t)=\sum_{k=0}^R \frac{t^k}{k!}$. Then, for any $R\in
\End(E)\otimes A_+$, and $H\in\herm(E)$, we have
$$
\Big|\!\Big| \e^{-(H+R)}\Big|\!\Big| \leq \e^{-\sm(H)}\Pcal(\|R\|).
$$
\end{lemm}

\begin{proof} Let $c=\sm(H)$.
Then $\|\e^{-uH}\|= \e^{-uc}$ for all $u\geq 0$. The term $I_k(H,R)$
of the Volterra expansion vanishes for $k> R$ since the term
$\e^{s_1H}R\cdots R\e^{s_{k+1}H}$ belongs to $\End(E)\otimes A_{k}$.
The norm of the term $I_k(H,R)$ is bounded by
$\frac{1}{k!}\e^{-c}\|R\|^k$. Summing up in $k$, we  obtain our
estimate.
\end{proof}

\bigskip

The preceding estimates hold if we work in the algebra
$\End(E)\otimes A$, where $E$ is a super-vector space and $A$ a
super-commutative algebra.

\end{document}